\documentclass{article}
\usepackage{amsthm,amsmath,amssymb,longtable,euler,latexsym}
\newcommand{\ds}{\displaystyle}
\newtheorem{theorem}{Theorem}[subsection]
\newtheorem{lemma}{Lemma}
\newtheorem{proposition}[theorem]{Proposition}

\theoremstyle{definition}

\theoremstyle{remark}

\newcommand{\tp}[1]{\ensuremath{\mathbf{#1}}}

\renewcommand{\strut}{{\rule{0pt}{2.5ex}}}
\renewcommand{\labelenumi}{\setlength{\labelwidth}{\leftmargin}
  \addtolength{\labelwidth}{-\labelsep} \hbox to
  \labelwidth{\theenumi.\hfill}} \numberwithin{equation}{subsection}
\title{Classification of Singular Fibres on Rational Elliptic Surfaces
  in Characteristic Three} 
\author{Tyler Jarvis, William E. Lang,\\
  Nansen Petrosyan, Gretchen Rimmasch,\\Julie Rogers, Erin~D.~Summers}
\date{\today}

\begin{document}

\maketitle
\section{Introduction}

\subsection{Background and Overview}
In~\cite{persson-fibreconfig} U. Persson classified all possible
configurations of singular fibres on rational elliptic surfaces over
the complex numbers.  He found $279$ different configurations.  R.
Miranda~\cite{miranda-singfibres} used a more combinatorial method to
redo Persson's list.  W. Lang~\cite{Lang-I} continued the
classification for characteristic two. Following the method of Lang,
we list all possible configurations of singular fibres on rational
elliptic surfaces in characteristic three.  In total, we find that $267$ distinct configurations exist.

Let $k$ be an algebraically closed field of characteristic three.
Let $f:X\to \mathbb{P}^1$ be a rational elliptic surface over $k$.
We will assume that the surface is relatively minimal and that $f$
has a section.  Such a surface has an expression in terms of a
Weierstrass equation
\begin{equation}y^{2}=x^{3}+ b_{2}x^{2}-b_{4}x + b_{6},
\end{equation} with $b_{2}$, $b_{4}$, $b_{6}\in k[t]$ and degree $b_{i}\leq
i$. The discriminant $\Delta$ for this surface is
\begin{equation}
 \Delta = -b_{2}^{2}(b_{2}b_{6}-b_{4}^{2})+b_{4}^{3}.
 \end{equation}
Assuming
that the surface has good reduction at infinity, the discriminant
$\Delta$ of the surface is a polynomial in $t$ of degree exactly
$12$. According to Tate's algorithm~\cite{Tate}, the surface has a
singular fibre of multiplicative type at a given point $t=p$ of
the base if and only if $\Delta (p)=0$ but $b_{2}(p)\neq0$, and a
singular fibre of additive type at $t=p$ if and only if
$\Delta(p)=0$ and $b_2(p)=0$.  There are four general
configurations of singular fibres that we will consider:
\begin{enumerate}
\item $b_2\neq 0$ and all the singular fibres are multiplicative.
\item $b_2\neq0$ and one of the singularities is additive, and any
  remaining singularities are multiplicative.
\item $b_2\neq0$ and two of the singularities are additive, and any
  remaining singularities are multiplicative.
\item $b_2\equiv 0$, in which case all the singularities are additive,
  and there can be more than two distinct additive singularities.  In
  this case the $j$-invariant is also identically zero.
\end{enumerate}

We will treat the first case separately and the other three
together. The structure of the paper will be as follows.  First,
some preliminaries relating to all surfaces will be discussed.  We
will then consider those surfaces which have only multiplicative
type singular fibres, followed by those surfaces with at least one
additive type singular fibre.  Finally we will give a summary of
the results.

\subsection{Acknowledgments}

We are grateful to Heidi Jarvis for her expert help with proofreading
and typesetting.  Research of the first author was supported in part
by NSF grant DMS-0105788.  The last four authors were supported in
part by graduate research assistantships from the Brigham Young University
mathematics department.

\section{Preliminaries}
\label{sec:pre} Let $f:X\to \mathbb{P}^1$ be a rational elliptic
surface over an algebraically closed field $k$ of characteristic
three.  We assign several invariants to each singular fibre, as given
in Table~\ref{tab:fibreinfo}.  First we list the Kodaira type $F$.
But, unlike in characteristic zero, each singularity may appear in
different ways, for example with different orders of vanishing of
$\Delta$, so we also list Lang's case number $L$ from the
classification in \cite{Lang-II} of isolated additive type singular
fibres on a rational elliptic surface in characteristic three.  Lang's
normal forms from that classification play a key role in this
project.  

The order of vanishing of $\Delta$ is denoted by $\delta$.  Associated
with each singular fibre is a lattice $\Lambda$ of rank $r$, generated
by the irreducible components of the singular fibre which do not meet
the zero section.  The lattice $\Lambda$ is a sublattice of the
N\`eron-Severi lattice of the entire surface, with discriminant $d$,
and is orthogonal to the lattices associate to the other fibres.

\begin{longtable}[c]{|c|c|c|c|c|c|}
\caption{Fibre Invariants \label{tab:fibreinfo}}
\\
\hline
\textbf{F}   & \textbf{L} & $\boldsymbol{\delta}$  & $\boldsymbol{r}$   & $\boldsymbol{\Lambda}$ & $\boldsymbol{d}$  \\
\hline \hline \endfirsthead
        \multicolumn{6}{l}{\small\slshape continued from previous page}\\
        \hline
\textbf{F}   & \textbf{L} & $\boldsymbol{\delta}$  & $\boldsymbol{r}$   & $\boldsymbol{\Lambda}$ & $\boldsymbol{d}$  \\
 \hline   \hline
\endhead
        \hline
        \multicolumn{6}{r}{\small\slshape continued on next  page}\\
        \endfoot
\endlastfoot
$\textrm{I}_n$ & & $n$ & $n-1$ & $A_{n-1}$ & $n$  \\
\hline
$\textrm{II}$ & 1A & 3 & 0 & $\{ 0 \} $ & 0 \\
\hline
$\textrm{II}$ & 1B & 4 & 0 & $\{ 0 \} $ & 0 \\
\hline
$\textrm{II}$ & 1C & 6 & 0 & $\{ 0 \} $ & 0 \\
\hline
$\textrm{II}$ & 1D & 7 & 0 & $\{ 0 \} $ & 0 \\
\hline
$\textrm{II}$ & 1E & 9 & 0 & $\{ 0 \} $ & 0 \\
\hline
$\textrm{II}$ & 1F & 12 & 0 & $\{ 0 \} $ & 0 \\
\hline
$\textrm{III}$ & 2 & 3 & 1 & $A_1 $ & 2 \\
\hline
$\textrm{IV}$ & 3A & 5 & 2 & $A_2 $ & 3 \\
\hline
$\textrm{IV}$ & 3B & 6 & 2 & $A_2 $ & 3 \\
\hline
$\textrm{IV}$ & 3C & 8 & 2 & $A_2 $ & 3 \\
\hline
$\textrm{IV}$ & 3D & 9 & 2 & $A_2 $ & 3 \\
\hline
$\textrm{IV}$ & 3E & 12 & 2 & $A_2 $ & 3 \\
\hline
$\textrm{I}_0^*$ & 4A & 6 & 4 & $D_4 $ & 4 \\
\hline
$\textrm{I}_0^*$ & 4B & 6 & 4 & $D_4 $ & 4 \\
\hline
$\textrm{I}_1^*$ & 5A & 7 & 5 & $D_5 $ & 4 \\
\hline
$\textrm{I}_2^*$ & 5B & 8 & 6 & $D_6 $ & 4 \\
\hline
$\textrm{I}_3^*$ & 5C & 9 & 7 & $D_7 $ & 4 \\
\hline
$\textrm{I}_4^*$ & 5D & 10 & 8 & $D_8 $ & 4 \\
\hline
$\textrm{IV}^*$ & 6A & 9 & 6 & $E_6 $ & 3 \\
\hline
$\textrm{IV}^*$ & 6B & 10 & 6 & $E_6 $ & 3 \\
\hline
$\textrm{IV}^*$ & 6C & 12 & 6 & $E_6 $ & 3 \\
\hline
$\textrm{III}^*$ & 7 & 9 & 7 & $E_7 $ & 2 \\
\hline
$\textrm{II}^*$ & 8A & 11 & 8 & $E_8 $ & 1 \\
\hline
$\textrm{II}^*$ & 8B & 12 & 8 & $E_8 $ & 1 \\
\hline
\end{longtable}
These invariants exclude certain configurations of
singular fibres, as they must satisfy the conditions given in the
following three lemmas.
\begin{lemma}[\cite{miranda-singfibres}]
\label{sigma-r}
On a rational elliptic surface, $\sum \delta_F = 12$ and $\sum r_F
\leq 8$.
\end{lemma}

\begin{lemma}[\cite{mirper}]
\label{prod-d}
If $\sum r_F = 8$, then $\prod d_F$ is a perfect square.
\end{lemma}

\begin{lemma}[\cite{shioda}]
\label{lattice} In order for a given configuration of singular
fibres to exist, the associated lattice $\oplus_F \Lambda_F$ must
embed in the root lattice $E_8$.
\end{lemma}

Oguiso and Shioda~\cite{oguiso-shioda} provide a list, originally due
to Dynkin~\cite{dynkin}, of lattices which embed in the root lattice
$E_8$.  For the convenience of the reader we include the complement of
that list here; namely, we list all direct sums of lattices of type
$A_n$, $D_n$ or $E_n$ of rank no more than eight which do not embed in
$E_8$.  Note that if the rank of such a lattice is $6$ or less, it
always embeds in $E_8$.

\begin{longtable}[c]{|l|l|}
\caption{Lattices which are direct sums of $A_n$'s, $D_n$'s,or $E_n$'s, which do not embed in $E_8$.
\label{tab:rootlattices}}
\\
\hline

\textbf{rank 7}   & \textbf{rank 8}    \\
\hline 
\endfirsthead
        \multicolumn{2}{l}{\small\slshape continued from previous page}\\
        \hline
\textbf{rank 7}   & \textbf{rank 8}   \\
    \hline
\endhead
        \hline
        \multicolumn{2}{r}{\small\slshape continued on next  page}\\
        \endfoot
\endlastfoot 
\hline
\strut $A_2 \oplus A^{\oplus 5}_{1}$ & $A_{2} \oplus A^{\oplus 6}_{1}$\\

 \strut          & $A_2^{\oplus 2} \oplus A_1^{\oplus 4}$\\

  \strut         & $A_3 \oplus A_1^{\oplus 5}$  \\
\hline
 \strut $A_2^{\oplus 2} \oplus A_1^{\oplus 3}$ & $A_3 \oplus A_2 \oplus A_1^{\oplus 3}$\\

 \strut          & $A_2^{\oplus 3} \oplus A_1^{\oplus 2}$  \\
\hline
 \strut  $A_3 \oplus A_2^{\oplus 2}$ & $A_3 \oplus A_2^{\oplus 2} \oplus A_1$ \\
\nopagebreak \strut     & $A_3^{\oplus 2} \oplus A_2$\\
\nopagebreak
\nopagebreak \strut  & $A_4 \oplus A_2^{\oplus 2}$\\*
\nopagebreak
\nopagebreak \strut  & $A_5 \oplus A_3$\\*
\nopagebreak
\nopagebreak \strut  & $A_6 \oplus A_2$\\
 \hline
 \strut  $A_4 \oplus A_1^{\oplus 3} $ & $A_4 \oplus A_1^{\oplus 4}$\\

 \strut  & $A_4 \oplus A_2 \oplus A_1^{\oplus 2}$\\

 \strut  & $A_4 \oplus A_3 \oplus A_1$\\
\nopagebreak
\nopagebreak \strut  & $A_5 \oplus A_1^{\oplus 3}$\\* 
\nopagebreak \strut  & $A_6 \oplus A_1^{\oplus 2}$\\*
\hline
\strut  $D_4 \oplus A_2 \oplus A_1$ & $D_4 \oplus A_2 \oplus A_1^{\oplus 2}$\\

 \strut  & $D_4 \oplus A_2^{\oplus 2}$\\

 \strut  & $D_4 \oplus A_3 \oplus A_1$\\

 \strut  & $D_5 \oplus A_2 \oplus A_1$\\

 \strut  & $D_6 \oplus A_2$\\

 \strut  & $D_7 \oplus A_1$\\
\hline
 \strut  & $D_5 \oplus A_1^{\oplus 3}$\\

 \strut  & $E_6 \oplus A_1^{\oplus 2}$\\
\hline
\end{longtable}

We also note that we have not used the J-map in our classifications at
all, although it was used in the classifications done in other
charactersitics.  We do not use it here because there were only a few cases
where the J-map provided simple results, and each of these cases could
be ruled out using another lemma which was required for more cases.

\section{Surfaces with Singular Fibres of Multiplicative Type}

We now consider the case where all singular fibres are of multiplicative type.

\subsection{General Information}

For convenience, when we write out a configuration of multiplicative
fibres, we will write $n_1\, n_2\, n_3 \dots$ instead of writing
$I_{n_1}I_{n_2}I_{n_3}\dots$.  Repeated terms will be denoted by
exponents, so $5\,3^2\,1$ denotes the configuration $I_5I_3I_3I_1$.

As mentioned before, a rational elliptic surface has a multiplicative
fibre over a point $t=p$ precisely if $\Delta$ has a root at $p$ but
$b_2$ does not.  Since $b_2$ is a non-zero polynomial of degree two,
we have two cases: either $b_2$ has distinct roots, or $b_2$ has a
repeated root.  In all cases where a given configuration exists, we
will give an example with distinct roots of $b_2$.  It is interesting
to note that there are 11 cases which exist when $b_2$ has distinct
roots but do not exist when $b_2$ has a repeated root; namely \tp{5\,3\,1^4},
\tp{8\,1^4}, \tp{5\,2^2 1^3}, \tp{5\,4\,1^3}, \tp{3^2 2^2 1^2}, \tp{6\,2^2 1^2}, \tp{5\,3\,2\,1^2}, \tp{8\,2\,1^2},
\tp{5^2 1^2}, \tp{4\,3\,2^2 1}, and \tp{4^2 2^2} .

We begin by reducing the discriminant using standard transformations
to a normal form, given in equations~(\ref{eq.*}) and (\ref{eq.**}).
In the case that $b_2$ has distinct roots, we use a fractional linear
transformation to move these roots to zero and infinity, and in the
case that $b_2$ has a double root, we move that unique root to
infinity.

\begin{proposition}
Consider a minimal rational elliptic surface over $\mathbb{P}^1_t$ 
with Weierstrauss equation
\begin{equation}\label{one}
y^2 = x^3 + b_{2}x^2 - b_{4}x + b_6,
\end{equation}
where the $b_{i}$ are polynomials in t of degree no more than $i$.  If all of the singular fibres are of multiplicative type, and if
the roots of $b_{2}$ are distinct, then up to automorphism of the
surface and the base, the discriminant of the surface is a degree
twelve polynomial of the form 
\begin{equation}\label{eq.*}
\Delta =
t^{12}+lt^{10}+t^3P_6(t)+nt^2+m,
\end{equation}
where $m^3=n^2l^3\neq 0$, and
where $P_6(t)$ is a polynomial of degree 6.

Conversely, given any $\Delta$ of the form (\ref{eq.*})
there is a rational elliptic surface with Weierstrass equation
(\ref{one}) and with roots of $b_2$ at $0$ and $\infty$.  If such
a surface exists with $\Delta$ of form (\ref{eq.*}), we may, up to
automorphism of the base, assume that one root of $\Delta$ occurs
at $t=1$.
\end{proposition}

\begin{proof} PGL(2) acts three-transitively on $\mathbb{P}^{1}$, thus we
can find a transformation that sends one root of $b_2$ to infinity and one root to zero.
Therefore, by a standard change of variables, we may assume $b_2 = t$.

 Let $$b_4=\xi_4t^4+\xi_3t^3+\xi_2t^2+\xi_1t+\xi_0,$$ and
$$b_6=\gamma_6t^6+\gamma_5t^5+\gamma_4t^4+\gamma_3t^3+\gamma_2t^2+\gamma_1t+\gamma_0,$$ where $\xi_i, \gamma_j
\in k$. By substituting $b_2, b_4, b_6$ into the equation of
$\Delta$ and simplifying we get
\begin{align*}
\Delta = &\xi_4^3 t^{12} + \xi_4^2t^{10}+ (-\gamma_6 - \xi_3
\xi_4)t^9 + (\xi_3^2 - \gamma_5 - \xi_2 \xi_4)t^8\\
& + (-\gamma_4 - \xi_1 \xi_4 - \xi_2 \xi_3)t^7
 + (-\gamma_3 + \xi_2^3 - \xi_0\xi_4 - \xi_1\xi_3 + \xi_2^2)t^6 \\
 &+ (-\gamma_2 - \xi_0 \xi_3 - \xi_1
\xi_2)t^5 + (-\gamma_1 - \xi_1^2 - \xi_0 \xi_2^2)t^4\\
 &+ (-\gamma_0 - \xi_0\xi_1 + \xi_1^3)t^3 + \xi_0^2 t^2 + \xi_0^3.
\end{align*}
So
$$\Delta = a^3t^{12} + a^2t^{10} + t^3P_6 + b^2t^2 + b^3,$$
for some $a, b \in k$, and $P_6 \in k[t]$ is a polynomial of degree
6 in t.  Since the surface has good reduction at infinity, we have
$a \neq 0$; so using the transformation
$x=a^{\frac{1}{2}}x'$ and $y=a^{\frac{3}{4}}y'$, we obtain
$$\Delta=(t^{12} +lt^{10}+t^3P_6+mt^2+n),
$$
where $$l=a^{-1},\quad m=b^2a^{-3}, \quad n=b^3a^{-3},$$ and
hence $m^3=n^2l^3$.

Conversely, for any $l, m,n \in k$ such that $m^3=n^2l^3 \neq 0$ and
$P_6 \in k[t]$, we may invert the previous steps; first, by using the
standard transformation $x\acute{}=a^{-\frac{1}{2}}x$ and
$y\acute{}=a^{-\frac{3}{4}}y$ where $a=l^{-1}$, and $b=nm^{-1}$,
giving $\Delta=a^3t^{12} + a^2t^{10} + t^3P_6 + b^2t^2 + b^3$.
Letting $$\xi_4 = a , \quad \xi_0 = b, \quad {\epsilon}_0 = {\epsilon}_2
=0, \quad {\epsilon}_1 = 1$$
and $\ds P_6(t) = \sum_{i=0}^6 \psi_it^i, $
with $\psi_i \in k$, we may choose any $\xi_1, \xi_2$ and $\xi_3 \in
k$ and solve for the coefficients of $b_6$:
$$\gamma_6 = -\xi_3 \xi_4 - \psi_6,\quad \gamma_5 = \xi_3^2 - \xi_2
\xi_4 - \psi_5,\quad \gamma_4 = -\xi_1 \xi_4 - \xi_2 \xi_3 - \psi_4,$$
$$\gamma_3 = \xi_2^3 - \xi_0 \xi_4 - \xi_1 \xi_3 + \xi_2^3 -
\psi_3,\quad \gamma_2 = -\xi_0 \xi_1 - \xi_1 \xi_2 - \psi_2,$$
$$\gamma_1 = -\xi_1^2 - \xi_0 \xi_2^3 - \psi_1, \quad \gamma_0 =
-\xi_0 \xi_1 + \xi_1^3 - \psi_0.$$

Substituting $b_2, b_4$, and $ b_6$ in $\Delta$, we obtain
$$\Delta = a^3 t^{12} + a^2 t^{10} + t^3 P_6 + b^2t^2 + b^3.$$
Finally, note that while PGL(2) is 3-transitive, we only moved two
points (the roots of $b_2$) to specific locations, so we may use
PGL(2) to move one root of $\Delta$ to $t=1$ as well. \end{proof}

A similar argument proves the next proposition.
\begin{proposition}
Consider a minimal rational elliptic surface over $\mathbb{P}^1_t$
with Weierstrauss equation
\begin{equation}\label{two}
y^2 = x^3 + b_{2}x^2 - b_{4}x + b_6,
\end{equation}
where $b_{i}$ are polynomials in t and degree $b_{i}\leq i$.  If
all of the singular fibres are of multiplicative type, and if $b_{2}$ has a double root, then up to automorphism of the surface
and the base, the discriminant $\Delta$ is a degree-twelve
polynomial of the form
\begin{equation}\label{eq.**}
\Delta = t^{12} + lt^{9} + mt^8- nt^7 + P_6,
\end{equation} where $l,m,n\in k$, such that $n^3=lm^3\neq 0$, and
$P_6(t)$ is a polynomial of degree less than or equal to 6.
Conversely, given any $\Delta$ of the form \ref{eq.**} there is a
rational elliptic surface with Weierstrass equation \ref{two} and
with roots of $b_2$ at $\infty$.  If such a surface exists, we may
assume one root of $\Delta$ occurs at $t=0$ and one occurs at
$t=1$.
\end{proposition}

To prove the existence of a rational elliptic surface with
all singular fibres of multiplicative type, we simply need to
exhibit a degree-twelve polynomial in the normal form 
with the zeros distributed in the desired way.
Conversely, if no such polynomials exist, no surface with the
prescribed fibre types exists.  This yields the following 
lemma for non-existence of configurations with fibres of high multiplicity.

\begin{lemma}
\label{mult3} A rational elliptic surface with purely
multiplicative fibres cannot have three (not necessarily distinct)
multiplicative fibres each of multiplicity at least three.  In particular, no purely multiplicative configuration may include $3^3$, $3^2 4$, $3\,4^2$, $4^3$, $3^2 5$, $3\,4\,5$, $3\,6$, $3\,7$, $3\,8$, $4\,6$, $5\,6$, $5\,7$, $6^2$, $9$, $10$, $11$, or $12$.
\end{lemma}

\begin{proof} Since we can
  move one root to $t=1$, the discriminant $\Delta$ of such a surface
  can be written in the form:
\begin{align} \label{expand}
\lefteqn{(t-1)^{3}(t-c)^{3} (t-d)^{3}(t^{3}+ e t^{2}+f t+ g)=} \\
& & t^{12}+e t^{11} + f t^{10} +(g-c^3-d^3-1)t^9
-(c^3e+d^3e+e)t^8 \nonumber \\
& & -(c^3f+d^3f-f)t^7
+(c^3-d^3g-(1+c^3)(g-d^3))t^6 \nonumber \\
& & +(c^3e+(1+c^3)d^3e)t^5  +
((1+c^3)d^3f+c^3f)t^4\nonumber \\
& & 
+((i+c^3)d^3g+c^3(g-d))t^3
 -c^3 d^3 e t^{2} -c^3 d^3 f t-c^3 d^3 g. \nonumber
\end{align}
If $b_2$ has distinct roots, $\Delta$ must be of the
form
$$t^{12}+lt^{10}+t^3P_6^{(t)}+mt^2 + nb^3,$$
 so the coefficient on the degree-$11$ term must vanish.  Thus $e=0$, and the polynomial is of the form
\begin{align*}
&t^{12}+ft^{10}+t^{3}P_{6}-c^{3}d^{3}f t-c^{3}d^3g.
\end{align*}
Since the coefficient on the degree-two term must be nonzero to
satisfy the form of Equation~(\ref{eq.*}), we have a contradiction.
Notice that no conditions were imposed on the distinctness of the
roots of the polynomial.

Suppose now that $b_2$ has a double root.  Since the coefficients
of the degree-$11$ and degree-$10$ terms must be zero, we have
$$e=f=0.$$
Thus we have
$$\Delta=t^{12}+(g-c^3-d^3)t^9+(c^3d^3-d^3g-c^3g)t^6+c^3d^3gt^3.$$ But
since the coefficients of the degree-$8$ and degree-$7$ terms must be
nonzero, no such polynomial satisfying Equation~(\ref{eq.**})
exists. Again, no conditions were imposed on the distinctness of
the roots. \end{proof}

\subsection{Solutions for Purely Multiplicative Singularities}

Below we list the 77 partitions of 12 using the notation of
\cite{miranda-singfibres} and \cite{Lang-I}.  
Next to each partition, we either list a
polynomial of the form given in Equation~(\ref{eq.*}) whose zeroes are distributed in the
way indicated by the partition, or we indicate that the surface
does not exist and give the letter of a lemma which proves its
non-existence (often there are several that apply).  As mentioned above, our results show that whenever a surface exists with purely multiplicative singularities and a double root of $b_2$, one also exists with distinct roots of $b_2$.  We know of no good proof of this fact other than an appeal to our list below.

In the following list $i$ denotes a root of $x^2+1=0$ in the field $k$.

\begin{enumerate}
\itemsep=.5ex
\item $\tp{1^{12}}$ \quad   $(t-1)(t^{11}+t^{10}-t^{9}+t-1)=t^{12}+t^{10}+t^9+t^2+t+1$

\item $\tp{2\, 1^{10}}$  \quad $(t-1)^2(t^{10}-t^{9}+t^{8}+t^{2}+t-1)= t^{12}+t^{10}+t^8+t^4-t^3+t^2-1$

\item $\tp{3\,1^9}$ \quad   $(t-1)^3(t^9+t^7+t^6+t^3-t^2+1)=t^{12}+t^{10}-t^7-t^5+t^2-1$

\item $\tp{2^2\,1^8}$ \quad  $(t^{2}+1)^{2}(t^{8}-t^{6}+t^{4}-t^{2}+1)=t^{12}+t^{10}+t^2+1$

\item $\tp{4\,1^8}$ \quad   $(t-1)^{4}(t^{8}+t^{7}-t^{6}-t^{3}-t^{2}+t+1)=t^{12}+t^{10}+t^7-t^6-t^5+t^4-t^3+t^2+1$

\item $\tp{3\,2\,1^7}$ \
$(t-1)^{3}(t+1)^{2}(t^{7}+t^{6}+t^{5}-t^{2}+t+1)=t^{12}+t^{10}-t^9+t^8+t^7-t^6+t^5+t^4-t^3+t^2-1$

\item $\tp{5\,1^7}$ \quad   $(t-1)^5(t^7-t^6+t^5+t^4-t^2+t-1)=t^{12}+t^{10}-t^8-t^7-t^6-t^3+t^2+1$

\item $\tp{2^3\,1^6}$ \quad$(t-1)^{2}(t^{2}+1)^{2}(t^6-t^5-t^4+at^3+t-1)$ \\
    $= t^{12}+t^{10}+(a+1)t^9+(a+1)t^8-at^6+(a+1)t^4+(a+1)t^3+t^2-1$\\
    where $a \neq 0,\pm 1$

\item $\tp{4\,2\,1^6}$ \quad   $(t-1)^4(t+1)^2(t^6-t^5-t^2+t-1)=t^{12}+t^{10}-t^9-t^7+t^6-t^5+t^2-1$

\item $\tp{3^2\,1^6}$ \quad  $(t-1)^3(t+1)^3(t^6+t^4-t^2+1)=t^{12}+t^{10}-t^8-t^4+t^2-1$

\item $\tp{6\,1^6}$ \quad  $(t-1)^6(t^6+t^4+t^2+1)=t^{12}+t^{10}+t^9+t^8+t^7-t^6+t^5+t^4+t^3+t^2+1$

\item $\tp{3\,2^2\,1^5}$ \quad$(t-1)^3(t^2+1)^2(t^5-t^3+1)=t^{12}+t^{10}-t^9-t^8-t^6-t^4-t^3+t^2-1$

\item $\tp{5\,2\,1^5}$ \quad$(t-1)^5(t+i + 1)^2(t^5+i t^4+t^3+i t^2+(1+ i )t+1)$ \\
    $= t^{12}+(i-1)t^{10}+i t^9+t^8-(i+1)t^7+(1+i)t^5-t^4+(i-1)t^3+(1-i)t^2+i $

\item $\tp{4\,3\,1^5}$ \quad$(t-1)^4(t+1)^3(t^5+t^4-t^3-t-1)=t^{12}+t^{10}+t^9-t^8-t^4-t^3+t^2-1$

\item $\tp{7\,1^5}$ \quad$(t-1)^7(t^5+t^4-t^3+t^2+dt+d)=t^{12}+t^{10}+(d-1)t^8+t^7-dt^6+(d-1)t^5+t^4-(d+1)t^3+(d-1)t^2-d$\\
    where $d^3-d^2-1 = 0$

\item $\tp{2^4\,1^4}$ \quad$(t-1)^2(t+1)^2(t^2+1)^2(t^4+t^2-1)=t^{12}+t^{10}+t^6+t^2-1$

\item $\tp{4\,2^2\,1^4}$ \quad$(t-1)^4(t^2+1)^2(t^4+t^3+t^2+t+1)=t^{12}-t^{10}-t^9+t^8+t^4-t^3-t^2+1$

\item $\tp{3^2\,2\,1^4}$ \quad$(t-1)^3(t+1)^3(t+i)^2(t^4+i t^3+t^2+i t-1)= t^{12}+t^{10}-i t^9-t^8-t^4+i t^3+t^2-1$

\item $\tp{6\,2\,1^4}$ \quad $(t-1)^6(t+1)^2(t^4+t^3+t^2+t+1)=t^{12}+t^{10}-t^9+t^8+t^7+t^5+t^4-t^3+t^2+1$

\item $\tp{5\,3\,1^4}$ \quad$(t-1)^5(t^3+1)(t^4-t^3+t^2-\alpha t
+\alpha)$\\
    $=t^{12}+t^{10}-i t^9+t^8+(i-1)t^6-t^4+i t^3-t^2-i$

\item $\tp{4^2\,1^4}$ \quad $(t^4+1)^2(t^4+t^2-1)=t^{12}+t^{10}+t^8-t^6-t^4+t^2-1$

\item $\tp{8\,1^4}$ \quad$(t-1)^8(t^4-t^3+t^2+t-1)=t^{12}+t^{10}-t^9+t^8+t^7+t^6+t^5+t^4+t^2-1$

\item $\tp{3\,2^3\,1^3}$ \quad$(t-1)^3(t+1)^2(t^2+ i t-1)^2(t^3+(1+i)t^2+(1-i)^2t+1-i)$\\
    $=t^{12}+t^{10}+t^9+t^8+(i-1)t^7+(1+i)t^6+(i-1)t^5-i t^4+(1+i)t^3-i t^2+i-1$

\item $\tp{5\,2^2\,1^3}$ \quad$(t-1)^5(t^2+i t+1)^2(t^3+(i-1)t^2+i(1-i)t-i$\\
    $ = t^{12}-i t^{10}+t^9-(1+i)t^8-t^7-i t^6+t^5+(i+1)t^4+t^3+i t^2+i$

\item $\tp{4\,3\,2\,1^3}$ \quad$(t-1)^4(t+1)^3(t+i)^2(t^3+(1+i)t^2+(1-i)t+1$\\
    $ = t^{12}+i t^{10}+ t^9 -i t^8- i t^4-t^3 +i t^2-1$

\item $\tp{7\,2\,1^3}$ \quad$(t-1)^7(t+i)^2(t^3+(i+1)t^2+c t-c(i+1))$\\
    $= t^{12}+(c-i-1)t^{10}+(i c+c)t^9+(1-c+i)t^8+(c-i-1)t^7+(1-c+i)t^5+(c-i-1)t^4-t^3+(1-c+i)t^2-c-i c$\\
    where $i c^5-c^3-(i+1)c^2-i+1 = 0$

\item $\tp{3^3\,1^3}$ \quad does not exist,~\ref{mult3}

\item $\tp{6\,3\,1^3}$\quad does not exist,~\ref{mult3}

\item $\tp{5\,4\,1^3}$ \quad$((t-1)^5(t+i)^4(t^3-(1+i)t^2+t-i-1)$\\
    $=t^{12}-(i-1)t^{10}+t^9-(i+1)t^8-i t^7+(1-i)t^6-i t^5+(1-i)t^4+(1-i)t^3+(1-i)t^2-1-i$

\item $\tp{9\,1^3}$\quad does not exist,~\ref{mult3}

\item $\tp{2^5\,1^2}$ \quad$(t^5-1)^2(t^2+1)=t^{12}+t^{10}+t^7+t^5+t^2+1$

\item $\tp{4\,2^3\,1^2}$ \quad$(t-1)^4(t^2+t+1)^2(t+1)^2(t^2+1)=t^{12}-t^{10}+t^9-t^8-t^7-t^6-t^5-t^4+t^3-t^2+1$

\item $\tp{3^2\,2^2\,1^2}$ \quad $(t-1)^3(t+1)^3(t^2+i)^2(t^2+i + 1)$\\
    $=t^{12}+t^{10}-i t^8+(1-i)t^6-t^4+i t^2+i+1$

\item $\tp{6\,2^2\,1^2}$ \quad $(t-1)^6(t^2+1)^2(t^2+i)$\\
    $= t^{12}+(i-1)t^{10}+t^9+(1-i)t^8+(i-1)t^7+(i+1)t^6+(1-i)t^5+(i-1)t^4+i t^3+(1-i)t^2+i$

\item $\tp{5\,3\,2\,1^2}$ \quad$(t-1)^5(t^3+i)(t+1)^2(t^2-1-i)$\\
    $=t^{12}-i t^{10}+(i-1)t^9-i t^8+(1+i)t^7+(i-1)t^6+(1+i)t^5-t^4-t^3-t^2+i-1$

\item $\tp{4^2\,2\,1^2}$ \quad$(t-1)^4(t-i)^4(t+1-i)^2(t^2-t+i)$\\
    $=t^{12}+t^{10}+i t^9+i t^8+(i-1)t^7-(1+i)t^6-(1+i)t^5-i t^4+t^3+t^2-1$

\item $\tp{8\,2\,1^2}$ \quad $(t-1)^8(t+1)^2(t^2+1)=t^{12}-t^{10}+t^9-t^8-t^7-t^6-t^5-t^4+t^3-t^2+1$

\item $\tp{4\,3^2\,1^2}$ \quad does not exist,~\ref{lattice}

\item $\tp{7\,3\,1^2}$ \quad does not exist,~\ref{prod-d}

\item $\tp{6\,4\,1^2}$ \quad does not exist,~\ref{lattice}

\item $\tp{5^2\,1^2}$ \quad$(t-1)^5(t+1)^5(t^2+i)$\\
    $=t^{12}+(1+i)t^{10}+(1+i)t^8+(i-1)t^6-(1+i)t^4-(1+i)t^2-i$

\item $\tp{10\,1^2}$ \quad does not exist,~\ref{sigma-r}

\item $\tp{3\,2^4\,1}$ \quad$(t-1)^3(t^4+t^3-t^2-t-1)^2(t+1)=t^{12}+t^{10}-t^8-t^7-t^6+t^3+t^2-1$

\item $\tp{5\,2^3\,1}$ \quad does not exist,~\ref{lattice}

\item $\tp{4\,3\,2^2\,1}$ \quad$(t-1)^4(t^3-i)(t^2+c)^2(t+1)$\\
    $=t^{12}-(1+c)t^{10}+(i-1)t^9+(c^2+c)t^8+(c-i+1-i c)t^7-(c^2+i)t^6+(i c^2-c^2+i c-c)t^5+(i c+i)t^4+(c^2-i c^2)t^3-(i c+ i c^2)t^2+i c^2$\\
    where $c^4-i c^3+c - i = 0$

\item $\tp{7\,2^2\,1}$ \quad does not exist,~\ref{prod-d}

\item $\tp{3^3\,2\,1}$\quad does not exist,~\ref{mult3}

\item $\tp{6\,3\,2\,1}$\quad does not exist,~\ref{mult3}

\item $\tp{5\,4\,2\,1}$\quad does not exist,~\ref{prod-d}

\item $\tp{9\,2\,1}$ \quad does not exist,~\ref{sigma-r}

\item $\tp{5\,3^2\,1}$\quad does not exist,~\ref{lattice}

\item $\tp{4^2\,3\,1}$ \quad does not exist,~\ref{lattice}

\item $\tp{8\,3\,1}$ \quad does not exist,~\ref{sigma-r}

\item $\tp{7\,4\,1}$ \quad does not exist,~\ref{sigma-r}

\item $\tp{6\,5\,1}$ \quad does not exist,~\ref{sigma-r}

\item $\tp{11\,1}$ \quad does not exist,~\ref{sigma-r}

\item $\tp{2^6}$ \quad$(t-1)^2(t^5+t^4+t+1)^2=t^{12}+t^{10}-t^6+t^2+1$

\item $\tp{4\,2^4}$ \quad$(t-1)^4(t^4-t^3+t^2-t+1)^2=t^{12}-t^{10}+t^9+t^7-t^6+t^5+t^3-t^2+1$

\item $\tp{3^2\,2^3}$ \quad does not exist,~\ref{lattice}

\item $\tp{6\,2^3}$ \quad does not exist,~\ref{lattice}

\item $\tp{5\,3\,2^2}$\quad does not exist,~\ref{lattice}

\item $\tp{4^2\,2^2}$ \quad$(t-1)^4(t+1)^4(t^2+1)^2=t^{12}+t^{10}-t^8+t^6-t^4+t^2+1$

\item $\tp{8\,2^2}$ \quad does not exist,~\ref{sigma-r}

\item $\tp{4\,3^2\,2}$ \quad does not exist,~\ref{prod-d}

\item $\tp{7\,3\,2}$  \quad does not exist,~\ref{sigma-r}

\item $\tp{6\,4\,2}$ \quad does not exist,~\ref{sigma-r}

\item $\tp{5^2\,2}$ \quad does not exist,~\ref{sigma-r}

\item $\tp{10\,2}$  \quad does not exist,~\ref{sigma-r}

\item $\tp{3^4}$  \quad does not exist,~\ref{mult3}

\item $\tp{6\,3^2}$  \quad does not exist,~\ref{sigma-r}

\item $\tp{5\,4\,3}$  \quad does not exist,~\ref{sigma-r}

\item $\tp{9\,3}$  \quad does not exist,~\ref{sigma-r}

\item $\tp{4^3}$ \quad  does not exist,~\ref{sigma-r}

\item $\tp{8\,4}$  \quad does not exist,~\ref{sigma-r}

\item $\tp{7\,5}$  \quad does not exist,~\ref{sigma-r}

\item $\tp{6^2}$  \quad does not exist,~\ref{sigma-r}

\item $\tp{12}$  \quad does not exist,~\ref{sigma-r}

\end{enumerate}

\section{At Least One Additive Fibre}

We now consider the case where the rational elliptic surface has at least one additive fibre.

\subsection{General Information}

In each of the cases below we will use a linear transformation 
to send the
worst singularity to $t=0$, using the fact that PGL(2) acts
three-transitively on $\mathbb{P}^1$.  We classify ``worst"
singularities using the following criteria:
\begin{enumerate}
\item Additive type singularities are worse than multiplicative
type singularities.
\item If comparing two additive type
singularities with different values of $\delta$ (where $\delta$
is the multiplicity of the root of $\Delta$ over which the
singularity happens), the singularity with the largest $\delta$ is
worse.
\item If two additive type singularities have the same
value of $\delta$, the singularity with the largest $r$ value is
worse.
\end{enumerate}

As in~\cite{Lang-II}, we may 
change coordinates so that the Weierstrass equation becomes
\begin{equation}\label{eq.normal}
y^2=x^3+tc_1x^2+tc_3x+tc_5,
\end{equation}
where the $c_i$ are polynomials in
$t$ of degree no more than $i$.  In our notation, $$b_2=tc_1,\quad b_4=-tc_3, \quad \mbox{and }
b_6=tc_5,$$ and $\Delta$ becomes
\begin{equation}\label{eq.delt-norm}
-\Delta = t^4c_1^2(c_1c_5-c_3^2)+t^3c_3^3.
\end{equation}
A list of all singular fibres of additive type that could potentially
appear on a rational elliptic surface in characteristic three was
given in \cite{Lang-II}, and a normal form over $t=0$ was given for
each one.  We will use the normal forms of \cite{Lang-II} for each of
the worst singularities.  For the reader's convenience, these normal
forms are listed at the beginning of each case treated in
Section~\ref{sec.AMsol}.

An easy application of Tate's algorithm (or a glance at
Equation~(\ref{eq.delt-norm})) shows that all additive fibres have
$\delta\geq3$.  The following lemma is also an easy result of Tate's
algorithm~\cite{Tate}.

\begin{lemma}
\label{b2} 
If $b_2\equiv 0$, then all singular fibres on the surface are of
additive type.  If $t^2|b_2$, where $b_2$ is not identically zero,
then there is only one additive singular fibre.
\end{lemma}

Since $PGL(2)$ is three-transitive, after sending the worst
singularity to $t=0$, we may still move two other points of the
base---either singular fibres or roots of $b_2$.

As in the purely multiplicative situation, we may
consider two cases: when $b_2$ has a double root at $t=0$, and when
$b_2$ has distinct roots.  If $b_2$ has distinct roots, we can move
the second root of $b_2$ to $t=\infty$.  Again, this gives some significant restrictions on $\Delta$.

\begin{proposition}\label{prp.add-delta}
  Consider a minimal rational elliptic surface over $\mathbb{P}^1_t$ with
  a singular fibre of additive type over $t=0$ and Weierstrass
  equation
\begin{equation}
y^2=x^3+tc_1x^2+tc_3x+tc_5,
\end{equation}
where the $c_i$ are polynomials in
$t$ of degree no more than $i$.  If the roots of $b_{2}$ are distinct, then up to automorphism of the
surface and the base, the discriminant of the surface is a polynomial
of degree 12 or less, of the form
\begin{equation}\label{eq.dr}
\Delta= -\alpha^3t^{12} + \alpha^2\beta^2 t^{10} + t^6P_3(t) + \beta^2 t^4 P_1(t) - \sigma^3 t^3, 
\end{equation}
where $\alpha$, $\beta$, and $\sigma$ are elements of $k$, and the
$P_i(t)$ are elements of $k[t]$ of degree no more than $i$.  In
particular, the $t^{11}$ term vanishes.

Similarly, if $b_2$ has a double root, then up to automorphism of the
surface and the base, the discriminant of the surface is a polynomial
of degree 12 or less, of the form
\begin{equation}\label{eq.rr}
\Delta=t^6P_6(t)  - \sigma t^3,
\end{equation}
where $\sigma$ is an element of $k$, and $P_6$ is an element of
$k[t]$ of degree no more than $6$.  In particular, the $t^4$ and $t^5$
terms vanish.

Finally, if $b_2$ is identically zero, then up to automorphism of the
surface and the base, the discriminant of the surface is a polynomial
of degree 12 or less, of the form
\begin{equation}\label{eq.b2=0}
\Delta = b_4^3 = -t^3(c_3)^3. 
\end{equation}
\end{proposition}

\begin{proof} We let
\begin{eqnarray}\label{eq.subs}
&b_2 = tc_1 &= t({\epsilon}_1 t+{\epsilon}_0)\\
&b_4 = -tc_3 &= t(\phi_3 t^3 + \phi_2t^2 + \phi_3t + \phi_0)\nonumber \\
&b_6 = tc_5 &= t(\gamma_5 t^5 + \gamma_4 t^4 + \cdots + \gamma_0),\nonumber
\end{eqnarray} 
with ${\epsilon}_i$, $\phi_i$ and $\gamma_i$ in $k$.  If $b_2$ has
distinct roots, then we may move the second root to infinity, which
amounts to setting  ${\epsilon}_1=0$. We substitute these into the
expression (\ref{eq.delt-norm}) for $\Delta$ to get 
\begin{align} 
\Delta&=-\phi_3^3t^{12} + {\epsilon}_0^2\phi_3^2t^{10}  - ({\epsilon}_0^2\phi_3\phi_2+{\epsilon}_0^3\gamma_5+\phi_2^3)t^9\\
&+(2{\epsilon}_0^3\gamma_4+{\epsilon}_0^2\phi_2^2+2{\epsilon}_0^2\phi_3\phi_1)t^8\nonumber \\
&- ({\epsilon}_0^3\gamma_3+{\epsilon}_0^2\phi_2\phi_1+{\epsilon}_0^2\phi_3\phi_0)t^7 \nonumber \\
& +({\epsilon}_0^2\phi_1^2+2\phi_1^3+2{\epsilon}_0^2\phi_2\phi_0+2{\epsilon}_0^3\gamma_2)t^6 \nonumber \\
&+ (2{\epsilon}_0^3\gamma_1+2{\epsilon}_0^2\phi_1\phi_0)t^5  +(2{\epsilon}_0^3\gamma_0+{\epsilon}_0^2\phi_0^2)t^4+2t^3\phi_0^3.\nonumber
\end{align}
This has the desired form.

Similarly, if $b_2$ has a double root, we may assume that ${\epsilon}_0=0$, and 
substituting into Equation~(\ref{eq.delt-norm}) gives
\begin{align} 
\Delta&=&({\epsilon}_1^2\phi_3^2+2{\epsilon}_1^3\gamma_5+2\phi_3^3)t^{12}+(2{\epsilon}_1^3\gamma_4+2{\epsilon}_1^2\phi_3\phi_2)t^{11}\\
&&+(2{\epsilon}_1^2\phi_3\phi_1+{\epsilon}_1^2\phi_2^2+2{\epsilon}_1^3\gamma_3)t^{10}  \nonumber \\
&&+(2{\epsilon}_1^2\phi_2\phi_1+2{\epsilon}_1^2\phi_3\phi_0+2{\epsilon}_1^3\gamma_2+2\phi_2^3)t^9
\nonumber \\
&&+(2{\epsilon}_1^2\phi_2\phi_0+2{\epsilon}_1^3\gamma_1+{\epsilon}_1^2\phi_1^2)t^8 \nonumber \\
&& +(2{\epsilon}_1^3\gamma_0+2{\epsilon}_1^2\phi_1\phi_0)t^7+({\epsilon}_1^2\phi_0^2+2\phi_1^3)t^6+2t^3\phi_0^3.\nonumber
\end{align}
This again has the desired form.

Finally, substituting $b_2=0$  into Equation~(\ref{eq.delt-norm}) gives the last form.
\end{proof}

An immediate corollary of the proposition is the following lemma.
\begin{lemma}
\label{add3} The configuration $\tp{X\,3^2\,2\,1}$ is not
possible, where $\tp{X}$ is any additive fibre.  Nor is any
specialization of this configuration possible, provided the $\tp{I_2}$ and
$\tp{I_1}$ fibres remain distinct from each other and both the order-three
multiplicative singularities remain distinct from the additive fibre $\tp{X}$. 
More precisely, the following configurations cannot
exist: \tp{X\,3^2\,2\,1}, \tp{X\,6\,2\,1}, \tp{X\,5\,3\,1},
\tp{X\,4\,3\,2}, \tp{X\,5\,4}, \tp{X\,3^2\,1}, \tp{X\,6\,1},
\tp{X\,3^2\,2}, \tp{X\,6\,2}, \tp{X\,4\,3}, \tp{X\,5\,3}.
\end{lemma}
\begin{proof}  Assume this configuration exists. 
  $X$ must have $\delta=3$, since $\sum \delta_i=12$, and we may
  assume $X$ occurs over $t=0$.  We may further assume no singular
  fibre occurs over infinity.  We let the first $I_3$ occur over
  $t=a$, the second $I_3$ over $t=b$, the $I_2$ over $t=c$, and the
  $I_1$ over $t=d$.
  
  The hypotheses of the lemma require $c\neq d$, and $ab\neq 0$ but
  do not otherwise require that the roots be distinct or non-vanishing.  We have
  $\eta t^3(t-a)^3(t-b)^3(t-c)^2(t-d)=\Delta$ for some $\eta\in k^*$.
We expand to get
\begin{multline*}
{t^3(t-a)^3(t-b)^3(t-c)^2(t-d)=}\\
t^{12}+(c+2d)t^{11}+(2cd+c^2)t^{10}\\
+(2b^3+2c^2d+2a^3)t^9+(2a^3c+2b^3c+a^3d+b^3d)t^8\\
+(2b^3c^2+a^3cd+b^3cd+2a^3c^2)t^7\\
+(a^3b^3+a^3c^2d+b^3c^2d)t^6+(2a^3b^3d+a^3b^3c)t^5\\
+(2a^3b^3cd+a^3b^3c^2)t^4+2t^3a^3b^3c^2d.
\end{multline*}
If $b_2$ has distinct roots, the $t^{11}$ term of $\Delta$ must
vanish, as in equation (\ref{eq.dr}), which means $c=d$, a
contradiction.

If $b_2$ has a repeated root, then the $t^5$ term of $\Delta$ must
vanish, as in equation (\ref{eq.rr}), which means $a^3b^3c=a^3b^3d$,
and thus either $c=d$ or $ab=0$, a contradiction.
\end{proof}

All of the configurations that do not exist fail to do so by one of
the previous lemmas, or by straightforward algebraic manipulations.
We list all that fail from algebraic manipulations in the following
lemma.

\begin{lemma}
\label{algman} The following configurations do not exist:
\begin{enumerate}
\item $\tp{II\, 3^2}$ and $\tp{II\, 6}$.
\item $\tp{IV\, 2^3}$ and $\tp{IV\, 4\, 2}$.
\item $\tp{IV\, 3\,2\,1}$ and $\tp{IV\, 5\,1}$.
\item $\tp{I_1^*\, 4\, 1}$.
\item $\tp{III\, 4^2\, 1}$ and $\tp{III\, 5\, 4}$.
\item $\tp{III\, III\, 3^2}$, $\tp{III\, III\, 6}$, $\tp{III\, II\,
    3^2}$, and $\tp{III\, II\, 6}$.
\item $\tp{III\, III\, 4\, 2}$, and $\tp{III\, II\, 4\, 2}$, where the
  $I_4$ and $I_2$ are distinct.
\end{enumerate}
\end{lemma}

\begin{proof}
\begin{enumerate}
\item $\tp{II\, 3^2}$ and $\tp{II\, 6}$.  Suppose this configuration
  exists.  Assume that one ${I_3}$ singularity occurs over $t=a$ and
  the other over $t=b$.  We allow the possibility of $a=b$.  The fibre
  of type $II$ must have $\delta=6$, which means it is case $1C$, which
  requires that ${\epsilon}_0=0$, $\phi_0=0$, and $\gamma_0 \neq 0$ (see
  Section~\ref{sec:1C}).  Making the substitutions of
  Equation~(\ref{eq.subs}) for $b_i$ in the discriminant, we have
\begin{align*}
  \alpha[(t-a)^3(t-b)^3] = & -\Delta/t^6 \\
  =&(2\phi_3^3+2{\epsilon}_1^3\gamma_5+{\epsilon}_1^2\phi_3^2)t^6
  +(2{\epsilon}_1^2\phi_3\phi_2+2{\epsilon}_1^3\gamma_4)t^5\\
  &+({\epsilon}_1^2\phi_2^2+2{\epsilon}_1^3\gamma_3+2{\epsilon}_1^2\phi_3\phi_1)t^4\\
  &+(2{\epsilon}_1^3\gamma_2+2{\epsilon}_1^2\phi_2\phi_1+2\phi_2^3)t^3\\
  &+({\epsilon}_1^2\phi_1^2+2{\epsilon}_1^3\gamma_1)t^2+
  2t{\epsilon}_1^3\gamma_0+2\phi_1^3,
\end{align*}
for some $\alpha \in k$.

Since the right-hand side of the above equation has a $t$ term, but
the left-hand side does not, the $t$ term must vanish.  But
${\epsilon}_1=0$ would force $b_2=0$, and thus all fibres would be
additive, and $\gamma_0=0$ is prohibited in case $1C$.

\item $\tp{IV\, 2^3}$ and $\tp{IV\, 4\, 2}$.  Suppose this
  configuration exists.  Since the type \tp{IV} fibre occurs with
  $\delta=6$, we are in case $3C$, with ${\epsilon}_1=0$,
  $\phi_0=0$, $\gamma_0=0$, but $\phi_1\neq0$.  Assume the
  multiplicative fibres occur over $a$, $b$, and $c$ in
  $\mathbb{A}^1-\{0\}$.  If all $a=b=c$, we have the configuration
  \tp{IV\, 6}, which \emph{is} possible (see Section~\ref{sec:3B}).
  If any two of these fibres are distinct, we may assume, without loss
  of generality, that $a=1$ and $b=-1$, which means we have
\begin{multline*}\alpha[t^6+ct^5+(1+c^2)t^4+ct^3+(1+c^2)t^2+ct+c^2] 
  =  -\Delta/t^6\\
  = (-{\epsilon}_1^3\gamma_5+{\epsilon}_1^2\phi_3^2-\phi_3^3)t^6
  +(2{\epsilon}_1^3\gamma_4-{\epsilon}_1^2\phi_3\phi_2)t^5
  +({\epsilon}_1^2\phi_2^2-{\epsilon}_1^2\phi_3\phi_1-{\epsilon}_1^3\gamma_3)t^4\\
  +(2\phi_2^3-{\epsilon}_1^2\phi_2\phi_1-{\epsilon}_1^3\gamma_2)t^3+(2{\epsilon}_1^3\gamma_1+{\epsilon}_1^2\phi_1^2)t^2-\phi_1^3.
\end{multline*}
Comparing low-degree coefficients shows that $c=0$, a contradiction.

\item $\tp{IV\, 3\,2\,1}$, $\tp{IV\, 5\,1}$.  Suppose this
  configuration exists.  Again the \tp{IV} falls in case $3C$, and we
  assume the $I_3$ fibre occurs over $a$, the $I_2$ over $b$, and the
  $I_1$ over $c$, with $abc\neq0$.  If the $I_2$ and $I_1$ fibres are
  distinct, we have $b\neq c$, and we may assume that $b=1$ and
  $c=-1$, which yields the following:
\begin{multline*}
  \alpha[{t}^{6}+2{t}^{5}+2{t}^{4}+ \left( 1+2{a}^{3} \right)
  {t}^{3}+{a}
  ^{3}{t}^{2}+{a}^{3}t+2{a}^{3]}=  -\Delta/t^6\\
  =   (-{\epsilon}_1^3\gamma_5+{\epsilon}_1^2\phi_3^2-\phi_3^3)t^6  
+(2{\epsilon}_1^3\gamma_4-{\epsilon}_1^2\phi_3\phi_2)t^5 
+({\epsilon}_1^2\phi_2^2-{\epsilon}_1^2\phi_3\phi_1-{\epsilon}_1^3\gamma_3)t^4\\
 +(2\phi_2^3-{\epsilon}_1^2\phi_2\phi_1-{\epsilon}_1^3\gamma_2)t^3 
+(2{\epsilon}_1^3\gamma_1+{\epsilon}_1^2\phi_1^2)t^2-\phi_1^3.
\end{multline*}
Comparing low-degree coefficients shows that $a=0$, a contradiction.

\item $\tp{I_1^* \,4\,1}$.  Suppose this configuration exists.  The
  $\tp{I_1^*}$ fibre falls into case $5A$, and thus has
  $\phi_0=\phi_1=\phi_2=0$ and $\gamma_0=\gamma_1=\gamma_2=0$, but
  ${\epsilon}_0\neq 0 $ and $\gamma_3\neq 0$.  In particular, $b_2$ has
  distinct roots. We may assume that the non-zero root of $b_2$ occurs
  at infinity, which implies that ${\epsilon}_1=0$, and the $I_4$ fibre
  occurs at $t=1$.  The remaining $I_1$ fibre occurs over some $a \in
  \mathbb{A}^1-\{0,1,\infty\}$.  Making these substitutions, we have
\begin{multline*}
\alpha[{t}^{5}- \left(1+a \right) {t}^{4}+{t}^{3}a-{t}^{2}+ \left( 1+a \right) t-a] = \Delta/t^7\\
=2{\phi_{{3}}}^{3}{t}^{5}+{\phi_{{3}}}^{2}{t}^{3}-{t}^{2}\gamma_{{
5}}+2t\gamma_{{4}}+2\gamma_{{3}}.
\end{multline*}
Since the coefficient of $t^4$ vanishes in $\Delta/t^7$, we have that $a=1$, which is a contradiction.

\item $\tp{III\, 4^2\, 1}$ and $\tp{III\, 5\, 4}$.
Suppose this configuration exists.  The type \tp{III} fibre falls into case $2$, which means that $\gamma_0=0$, but $\phi_0\neq0$ and ${\epsilon}_0\neq0$.   Let the $I_4$'s happen at
$t=a$ and $t=b$, and let the $I_1$ happen at $t=c$.  We require that $a\neq b$, but we could have $a=c$ or $b=c$. We have
\begin{equation}\label{eq:ddd}
\Delta =\alpha(t^3(t-a)^4(t-b)^4(t-c))
\end{equation}   
for some $\alpha \neq 0$.  If $b_2$ has two distinct roots, then by
Proposition~\ref{prp.add-delta} we may assume that ${\epsilon}_1=0$ and
$\Delta$ has no $t^{11}$ term.  This means $c=-a-b$.  Furthermore, we
may assume that $a=1$.  Making these substitutions into $\Delta$ shows
that the leading term vanishes if $\phi_3$ does, which would put an
additional additive singularity at infinity. Thus $\phi_3$ is
non-zero.  Equating like terms in Equation~(\ref{eq:ddd}), we obtain
the following:
\begin{align*}
{\frac {{\phi_{{0}}}^{3}}{{\phi_{{3}}}^{3}}}&={b}^{5}+{b}^{4}
  & \mbox{(from
$t^{3}$)}\\
-{\frac {{{\epsilon}_{{0}}}^{2}{\phi_{{0}}}^{2}}{{\phi_{{3}}}^{3}}}&=-
{b}^{4}-{b}^{3}-{b}^{5} &\mbox{(from $t^4$)} \\
-{\frac {{{\epsilon}_{{0}}}^{2}}{\phi_{{3}}}}&=-1-b-{b}^{2}
&\mbox{(from $t^{10}$).}
\end{align*}
A straightforward manipulation of these constraints
gives $b=1$, which is a contradiction, since $a=1$ and $a$ and $b$ are distinct.  

If $b_2$ has a repeated root, by Proposition~\ref{prp.add-delta} the $t^5$ and $t^4$ terms must vanish, and we may assume that $a=1$ and $b=-1$.  Substituting into Equation~(\ref{eq:ddd}), and setting like terms equal, we get that $c=0$, which is a contradiction.

\item $\tp{III\, III\, 3^2}$, $\tp{III\, III\, 6}$, $\tp{III\, II\, 3^2}$, and $\tp{III\, II\, 6}$.

Suppose this configuration is possible. We may assume that the type \tp{III} fibre occurs over  $t=0$, which implies that $\phi_0\neq 0$, but $\gamma_0=0$.  
We further assume that the other additive fibre occurs over $t=\infty$, so that $b_2$ has distinct roots, and we may take ${\epsilon}_0=1$ and ${\epsilon}_1=0$.
Furthermore, suppose the first $I_3$ happens at $t=a$ for some $a\neq 0, \infty$, and that
the second $I_3$ happens at $t=b\neq0, \infty$. Making these substitutions into $\Delta$ shows that the coefficient of $t^4$ in $\Delta$ is $\phi_0^2$. But matching coefficients in the equality  
$$\Delta=\alpha t^3(t-a)^3(t-b)^3=t^9-(a^3+b^3)t^6-t^3a^3b^3,$$
with $\alpha \in k^*$, shows that $\phi_0=0$, which is a contradiction.

\item $\tp{III\, III\, 4\, 2}$ and $\tp{III\, II\, 4\, 2}$. 
Suppose this configuration is possible. As in the previous case, we may assume that $\gamma_0={\epsilon}_1=0$ and ${\epsilon}_0=1$.  Since the $t^8$ term vanishes in the expansion of $\Delta$, matching coefficients in the equality $$\Delta = \alpha(t-a)^4(t-b)^2$$  shows that $a=b$, which is a contradiction.

\end{enumerate}
\end{proof}

\subsection{Solutions for Additive and Multiplicative Singularities}\label{sec.AMsol}

We list Lang's type \cite{Lang-II} for the worst singularity, its main
invariants, Weierstrass form, and $-\Delta$.  For each subsequent
configuration, we list the Weierstrass form and $-\Delta$ if the
desired rational elliptic surface exists, or the lemma number
proving impossibility if it does not exist.  In the following
examples, $i$ denotes a root of $x^2+1=0$ in $k$.

\subsubsection{Case 1A}
\label{subsubsec:1a}

$$\delta =3,\quad r=0\quad \mbox{ Type \tp{II}}$$
\begin{equation}
 y^2=x^3+tc_1x^2+tc_3x+tc_5, \quad t\nmid c_5, t\nmid c_3
 \end{equation}
$$-\Delta=t^2c_1^2(t^2c_1c_5-t^2c_3^2)+t^3c_3^3$$

\begin{enumerate}
\itemsep=.5ex
\item \tp{II\, 9} \quad $y^2 =x^3+t^2x^2+t(t+1)x+t(t+2)$\\
            $-\Delta =2t^3$

\item \tp{II\,  8 \, 1   } \quad does not exist,~\ref{add3}

\item \tp{II\,  7 \, 2    } \quad does not exist,~\ref{add3}

\item \tp{II\,  7 \, 1^2   } \quad $y^2 =x^3+t(t+i)x^2+t(2t+i)x+t^2$ \\
$-\Delta =-2it^3+t^4-2it^5 =-t^3(2+it+2t^2)i$

\item \tp{II\,  6 \, 3    } \quad $y^2 =x^3+t^2x^2+t(2t+2)x+t(2+t)$ \\
$-\Delta =t^9+2t^6=t^6(t^3+2)$

\item \tp{II\, 6\, 2 \, 1   } \quad does not exist,~\ref{add3}

\item \tp{II\,  6 \, 1^3  } \quad $y^2 =x^3+t(t+1)x^2+t(t+2+i)x+t(t+2i)$ \\
$-\Delta =(-1-i)t^6+it^5-it^4+(1+i)t^3=-t^3(t-1)((1+i)t^2+t+(1+i))$

\item \tp{II\,  5 \, 4    } \quad does not exist,~\ref{add3}

\item \tp{II\, 5\, 3 \, 1    } \quad does not exist,~\ref{add3}

\item \tp{II\,  5 \, 2^2   } \quad $y^2=x^3+t{\epsilon}_0x^2+t(2t^3+2+i)x+\frac{t}{{\epsilon}_0}((2+2i)t^5+2it^4+(2+i)t^3+t^2+(2+2i)t+2i)$\\
$-\Delta =t^{12}+it^{10}+(-1+i)t^9+2t^8+2it^6+(-1+i)t^5+t^4+(1+i)t^3=t^3(t+(1+i))^5(t-1)^2(t-i)^2$, where ${\epsilon}_0^2=i$

\item \tp{II\, 5\, 2 \, 1^2  } \quad $y^2 =x^3+tx^2+t(2t^3+i)x+t(it^5+t^4+2it^3+t^2+2it+2)$ \\
$-\Delta =t^{12}+t^{10}+2it^9+2t^8+2it^7+2t^6+it^5+it^3=t^3(t+i)^5(t-i)^2(t-1)(t+1)$

\item \tp{II\,  5 \, 1^4   } \quad $y^2=x^3+t{\epsilon}_0x^2+t(2t^3+i)x+\frac{t}{{\epsilon}_0}((1-i)t^5+(1+i)t^4+t^3+(1+i)t^2+it-i)$\\
$-\Delta =t^{12}+(1+i)t^{10}+t^9+it^8+t^7+it^6+(1+2i)t^5+t^4+it^3=t^3(t-(1-i))^5(t-i)(t+1)(t-1)(t-(1+i))$

\item \tp{II\, 4^2 \, 1 } \quad $y^2
  =x^3+t(t+{\epsilon})x^2+tx+t\left(\frac{-t}{{\epsilon}-1}+
    \frac{{\epsilon}}{{\epsilon}-1}\right)$ \\
  $-\Delta=\frac{t^8}{{\epsilon}-1}-\frac{{\epsilon} t^7}{{\epsilon}-1}
  +t^6-\frac{t^5}{{\epsilon}-1}+\frac{{\epsilon} t^4}{{\epsilon}-1}
  +\frac{(1-{\epsilon})t^3}{{\epsilon}-1}
  =t^3(t-1)^4\left(2+\frac{1}{{\epsilon}-1}
    -\frac{{\epsilon}}{{\epsilon}-1}-1+{\epsilon}^2\right)t$, where
  ${\epsilon}^3={\epsilon}^2-{\epsilon}-1$

\item \tp{II\, 4\, 3 \, 2    } \quad does not exist,~\ref{add3}

\item \tp{II\, 4\, 3 \, 1^2   } \quad $y^2 =x^3+t(t+1)x^2+tx+t(t+2)$ \\
$-\Delta =2t^8+t^7+t^6+t^5+2t^4+2t^3 =-t^3(t-1)^3(2-t+t^2)$

\item \tp{II\, 4\, 2^2 \, 1    } \quad $y^2 =x^3+tx^2i+t(2t^3+2i)x+t(2t^5+it^4-it^2+t+i)$ \\
$-\Delta =t^{12}+2t^{10}-it^9+2t^8+it^7+t^6+it^5+2it^3 =t^3(t+i)^4(t-1)^2(t+1)^2(t-i)$

\item \tp{II\, 4\, 2 \, 1^3   } \quad $y^2=x^3+t{\epsilon} x^2+t(2t^3+1)x+\frac{t}{{\epsilon}}(-t^5-(1+i)t^4-t^3+(1-i)t^2+t-1)$\\
$-\Delta=t^{12}+(2-i)t^{10}+(2+2i)t^9+it^8+(1+i)t^7+2t^6+(1+i)t^5+(1+i)t^4+2t^3=t^3(t+(1+i))^4(t+i)^2(t-(1+i))(t-(1-i))(t+1)$, where ${\epsilon}^2=2-i$

\item \tp{II\,  4 \, 1^5   } \quad $y^2=x^3+t{\epsilon} x^2+t(2t^3+1+i)x+\frac{t}{{\epsilon}}((1-2i)t^5-it^4+t^2+(2-i)t-2i)$\\
$-\Delta=t^{12}+it^{10}+(1-i)t^9+2t^8+(2+i)t^7-it^6+(2+i)t^5+2t^4+(2+i)t^3=t^3(t-(2-i))^4(t-1)(t+i)(t-i)(t+1)(t-(1+i))$, where ${\epsilon}^2=i$

\item \tp{II\, 3^3   } \quad $y^2 =x^3+t^2x^2+t(t+1)x+t(2+2t^2+t)$ \\
$-\Delta =t^9+2t^3=t^3(t-1)^3(t+1)^3$

\item \tp{II\, 3^2\, 2 \, 1    } \quad does not exist,~\ref{add3}

\item \tp{II\, 3^2 \, 1^3   } \quad $y^2 =x^3+t(t+1)x^2+t(2t+1)x+t(2+2t^2+2t)$ \\
$-\Delta =t^9+2t^8+t^7+t^5+2t^4+2t^3=-t^3(t-1)^3(-t^3+t^2-t-1)$

\item \tp{II\,  3 \, 2^3  } \quad $y^2 =x^3+t{\epsilon} x^2+t(2t^3+i)x+\frac{t}{{\epsilon}}((1-i)t^5-it^4-2it^2-it-i)$ \\
$-\Delta =t^{12}+it^{10}+(2+2i)t^9+2t^8+2t^7+t^6+2t^5+(2-i)t^4-2it^3 =t^3(t+1)^3(t-(2-i))^2(t-1)^2(t-i)^2$, where ${\epsilon}^2=i$

\item \tp{II\, 3\, 2^2 \, 1^2   } \quad $y^2 =x^3+t(t+1)x^2+t(t+2)x+t((2+i)t^2+(2+i)t+2+i)$ \\
$-\Delta =(1-i)t^9-(1+i)t^8+(1-i)t^7+(1-i)t^6+(1-i)t^5-(1+i)t^4+t^3=i\left(t-\frac{1}{1-i}\right)^2(2t^2+2it^2+1)t^3(t-1)^2$

\item \tp{II\, 3\, 2 \, 1^4   } \quad $y^2 =x^3+t^2x^2+t(t+1)x+t(1+t^2)$ \\
$-\Delta =2t^9+t^8+t^7+2t^3=t^3(t-1)^2(-t^4-t^3+t-1)$

\item \tp{II\,  3 \, 1^6   } \quad $y^2 =x^3+(t^2+it)x^2+(it^2+2t+t^2)x+2t^3+2t^2-it^2+it$ \\
$-\Delta =t^9+t^8+t^6+t^4+t^3=t^3(t^6+t^5+t^3+t+1)$

\item \tp{II\, 2^4 \, 1    } \quad $y^2 =x^3+t^2x^2+2tx+t(2t^5+2t^4+1+t^2+2t)$ \\
$-\Delta =t^{12}+t^{11}+2t^9+t^8+2t^7+t^6+t^3=t^3(t-1)^2(t^3+t^2+t+2)^2(t-2)$

\item \tp{II\, 2^3 \, 1^3  } \quad $y^2 =x^3+itx^2+t(2t^3+2+i)x+t((1+2i)t^5+2+(2+i)t^3+2t^2+2it)$ \\
$-\Delta =t^{12}+2t^{10}+(1+i)t^9+it^7-it^6+t^5+it^4+(1+i)t^3=t^3(t-(1+i))^2(t+(1+i))^2(t+1)^2(t-1)(t-i)(t-(1-i))$

\item \tp{II\, 2^2 \, 1^5  } \quad $y^2 =x^3+t{{\epsilon}}x^2+t(2t^3+2i)x+\frac{t}{{\epsilon}}((1-i)t^5-it^4+t^3+t^2+(-1+i)t+(1-i))$ \\
$-\Delta =t^{12}+(1+i)t^{10}+t^9+(2+i)t^8+it^7+(2+2i)t^6+2t^5+2it^4+2it^3t^3(t-(1+i))^2(t-1)^2(t+1)(t+i)(t-i)(t+(1+i))(t-(1-i))$, where ${\epsilon}^2=1+i$

\item \tp{II\,  2 \, 1^7  } \quad $y^2 =x^3+t{{\epsilon}}x^2+t(2t^3+2+i)x+\frac{t}{{\epsilon}}((-1+2i)t^5+it^4+(-2+2i)t^3+t^2+(-2+i)t+i)$ \\
$-\Delta=t^{12}-it^{10}+(4+2i)t^9+2t^8+(5+5i)t^7-2it^6+(2+4i)t^5+(4+4i)t^3 =t^3(t-(-1-i))^2(t-1)(t+1)(t-i)(t-(-1-i))(t+i)(t-(-1+i))(t-(1-i))$, where ${\epsilon}^2=-i$

\item \tp{II\, 1^9  } \quad $y^2 =x^3+t(t+1)x^2+2tx+t(t^4+2t^3+1)$ \\
$-\Delta =2t^{11}+t^{10}+2t^8+t^6+2t^5+t^3=-t^3(t-1)(1+t+t^3+t^4+t^7)$

\item \tp{II\, II\, II\, II  } \quad 
$y^2 = x^3-t(t^2-1)x+ t(t^2-1)(t^2-t-1)$\\
$-\Delta =t^9-t^3 = (t(t^2-1))^3$

\item \tp{II\, II\, 6  } \quad $y^2 =x^3+tx^2+t(2t^2+t+2)x+t(t^4+t^3+t^2+t+1)$ \\
$-\Delta =t^9+t^6+t^3=t^3(t-1)^6$

\item \tp{II\, II\,  5 \, 1   } \quad $y^2 =x^3+tx^2i+t(t^2+t+2)x+t^5i$ \\
$-\Delta =2t^9+t^8+t^7+2t^5+2t^4+t^3=-t^3(t-1)^5(t+1)$

\item \tp{II\, II\,  4 \, 2   } \quad $y^2 =x^3+tx^2i+t(2t^2+t+2)x+t(t^4i+t^2i+i)$ \\
$-\Delta =t^9+t^8+2t^7+t^6+2t^5+t^4+t^3=t^3(t-1)^4(t+1)^2$

\item \tp{II\, II\,  4 \, 1^2  } \quad $y^2 =x^3+tx^2+t(2t^2+t+2)x+t(2t^4+t^2+2)$ \\
$-\Delta =t^9+2t^8+t^7+t^6+t^5+2t^4+t^3=t^3(t-1)^4(t-\frac{1}{i})(t+\frac{1}{i})$

\item \tp{II\, II\, 3^2  } \quad $y^2 =x^3+tx^2+t(t^2+t+2)x+t(t^4+2t^3+t^2+1)$ \\
$-\Delta =2t^9+t^3=-t^3(t-1)^3(t+1)^3$

\item \tp{II\, II\, 3\, 2 \, 1   } \quad $y^2 =x^3+tx^2+t(t^2+t+2)x+t((2-i)t^4-t^3i+(1+i)t^2-ti+2i)$ \\
$-\Delta =2t^9+(2+i)t^8+(2+i)t^7+-it^6+(1+i)t^5+(1+i)t^4+t^3 =-t^3(t-\frac{1}{i})^3(t-1)^2(t+\frac{1}{i})$

\item \tp{II\, II\,  3 \, 1^3  } \quad $y^2 =x^3+tx^2+t(t^2+1)x+t(t^4+2t^3+2)$ \\
$-\Delta =2t^9+t^7+2t^6+2t^4+2t^3=-t^3(t-1)^3(t^3+2t+2)$

\item \tp{II\, II\, 2^3  } \quad $y^2 =x^3+tx^2i+t(2t^2+t+2)x+t(t^4i-t^3i-ti+i)$\\
$-\Delta =t^9+t^8+2t^6+t^4+t^3=t^3(t-1)^2(t-\frac{1}{i})^2(t+\frac{1}{i})^2$

\item \tp{II\, II\, 2^2 \, 1^2  } \quad $y^2 =x^3+tx^2+t(2t^2+t+2)x+t(t^4+2t^3+2t^2+2t+1)$ \\
$-\Delta =t^9+2t^7+2t^5+t^3=t^3(t-1)^2(t+1)^2(t-\frac{1}{i})(t+\frac{1}{i})$

\item \tp{II\, II\,  2 \, 1^4  } \quad $y^2 =x^3+tx^2+t(2t^2+t+2)x+t(t^4+t^3+2t^2+2t+2)$ \\
$-\Delta =t^9+2t^5+2t^4+t^3=t^3(t-1)^2(t^4-t^3+t+1)$

\item \tp{II\, II\, 1^6  } \quad $y^2 =x^3+tx^2+t(t^2+t+2)x+t(t^4+t^3+t^2+2t+1)$ \\
$-\Delta =2t^9+2t^7+2t^5+t^3=-t^3(t-1)^2(t+1)(t^4+1)$

\end{enumerate}

\subsubsection{Case 1B}
$$\delta=4,\quad r=0\quad \mbox{ Type \tp{II}}$$
\begin{equation}
y^2=x^3+tc_1x^2+t^2c_2x+tc_5, \quad t\nmid c_5, t\nmid c_1
\end{equation}
$$-\Delta=t^2c_1^2(t^2c_1c_5-t^4c_2^2)+t^6c_2^3$$

\begin{enumerate}
\itemsep=.5ex

\item \tp{II\, 8  } \quad does not exist,~\ref{add3}

\item \tp{II\,  7 \, 1   } \quad $y^2 =x^3+t(t+1)x^2+t^2x+t(2+t)$ \\
$-\Delta =2t^5+t^4=-t^4(t-1)$

\item \tp{II\,  6 \, 2   } \quad does not exist,~\ref{add3}

\item \tp{II\,  6 \, 1^2  } \quad $y^2 =x^3+t(t+1)x^2+2t^2x+t(2+t)$ \\
$-\Delta =2t^6+2t^5+t^4=t^4(-t^2-t+1)$

\item \tp{II\,  5 \, 3   } \quad does not exist,~\ref{add3}

\item \tp{II\, 5\, 2 \, 1    } \quad $y^2 =x^3+t(t+i)x^2+(1-i)t^2x+t(2+ti)$ \\
$-\Delta =2t^7+(-1+i)t^6+2t^5-it^4=t^4(-i-t-(1-i)t^2-t^3)$

\item \tp{II\,  5 \, 1^3   } \quad $y^2 =x^3+t(t+1)x^2+t^2x+t(t+1)$ \\
$-\Delta =t^7+2t^5+2t^4=t^4(2-t+t^3)$

\item \tp{II\, 4^2   } \quad $y^2 =x^3+t(t+1)x^2+t^2x+t(1+2t)$ \\
$-\Delta =2t^8+t^7+t^5+2t^4=-t^4(t-1)^4$

\item \tp{II\, 4\, 3 \, 1    } \quad $y^2 =x^3+t(t+1)x^2+t^2x+t((2-i)t+1+i)$ \\
$-\Delta =(-1+i)t^8-(2+i)t^7+(-2+i)t^5-(1+i)t^4 =-(1-i)t^4(t-1)\left(t-\frac{1}{i}\right)^3$

\item \tp{II\,  4 \, 2^2  } \quad $y^2 =x^3+t(t+i)x^2+t^2x-it$ \\
$-\Delta =t^8+t^6+t^4=t^4(t-1)^2(t+1)^2$

\item \tp{II\, 4\, 2 \, 1^2  } \quad $y^2 =x^3+t(t+1)x^2+t^2xi+t(ti+1+i)$ \\
$-\Delta =(2+2i)t^8-it^7+(2+i)t^6-it^5+(2+2i)t^4=t^4(t-1)^2((2+2i)t^2+t+(2+2i))$

\item \tp{II\,  4 \, 1^4   } \quad $y^2 =x^3+t(t+1)x^2+2t^2x+t$ \\
$-\Delta =t^8+t^7+2t^6+2t^4=-t^4(2-t+t^2-t^3)(t-1)$

\item \tp{II\, 3^2 \, 2   } \quad  does not exist,~\ref{add3}

\item \tp{II\, 3^2 \, 1^2   } \quad $y^2 =x^3+t(t+1)x^2+t^2xi+t(2+(1-i)t^2+t)$ \\
$-\Delta =(-1+i)t^9+t^8+2t^7+(1-i)t^6+2t^5+t^4=-t^4(t-1)^3((1-i)t^2-t+1)$

\item \tp{II\, 3\, 2^2 \, 1    } \quad $y^2 =x^3+t(t+i)x^2+2t^2x+t(1-(1+i)t^2+(1+i)t)$ \\
$-\Delta =(1+i)t^9+2it^8-(1+i)t^7+(1-i)t^6+(2+i)t^5+t^4i=it^4(t-1)^2\left(t-\frac{1}{1-I}\right)^2((1-i)t+i)$

\item \tp{II\, 3\, 2 \, 1^3  } \quad $y^2 =x^3+t(t+1)x^2+t^2x+t(2+t^2)$ \\
$-\Delta =2t^9+t^8+2t^6+t^4 =t^4(t-1)^2(1-t-t^2-t^3)$

\item \tp{II\,  3 \, 1^5   } \quad $y^2 =x^3+t(t+i)x^2+t^2x+t(1+(2-i)t^2)$ \\
$-\Delta =(-2+i)t^9+t^8+(2-i)t^7-(1+i)t^6+t^4i=-t^4(t-1)(t+1)(i-t^2+(2-i)t^3)$

\item \tp{II\, 2^4  } \quad $y^2 =x^3+t(t+{\epsilon})x^2+t^2x+t(-it^3+t^2({\epsilon}+{{\epsilon}}i)+(1-i)t+2i{\epsilon})$ \\
$-\Delta=t^{10}i+(2{\epsilon}-i{\epsilon})t^9+it^8+(i{\epsilon}+{\epsilon})t^7+2it^6+(2{\epsilon}-i{\epsilon})t^5+2it^4=it^4\left(t-\frac{1}{i}\right)^2\left(t+\frac{1}{i}\right)^2\left(t-\frac{{\epsilon}}{1+i}\right)^2$, where ${\epsilon}^4=-1$

\item \tp{II\, 2^3 \, 1^2  } \quad $y^2 =x^3+t(t+1)x^2+t^2x+t(2t^3+(1+i)t^2+ti+i)$ \\
$-\Delta =t^{10}+(-1+2i)t^9-(2+i)t^8-it^7+(2+2i)t^6+2it^5-it^4 =-t^4(t-1)^2\left(t-\frac{1}{i}\right)^2(-t^2-t-i)$

\item \tp{II\, 2^2 \, 1^4  } \quad $y^2 =x^3+t(t+1)x^2+2t^2x+t(t^3+2t^2+1) $\\
$-\Delta =2t^{10}+t^9+t^8+2t^4=t^4(t-1)^2(2+t-t^3-t^4)$

\item \tp{II\,  2 \, 1^6  } \quad $y^2 =x^3+t(t+1)x^2+t^2x+t(2+2t^3+2t^2)$ \\
$-\Delta =t^{10}+t^9+t^8+t^7+t^6+t^4=-t^4(t-1)\left(t-\frac{1}{i}\right)\left(t+\frac{1}{i}\right)(1+t+t^2-t^3)$

\item \tp{II\, 1^8   } \quad $y^2 =x^3+t(t+1)x^2+t^2x+t(t^4i-it^3+(1+i)t^2+(1-i)t+1)$ \\
$-\Delta =-it^{11}+t^{10}i-(1+i)t^9+(-2+i)t^7+(2-i)t^6+(-1+i)t^5+2t^4=-t^4(t-1)(t+i)(t-i)(2+(1+i)t+t^2+it^4)$

\item \tp{II\, II\, 5   } \quad $y^2 =x^3+tx^2+t^3x+t(2t^4+t^3+2t^2+2t^2+2)$ \\
$-\Delta =2t^9+2t^8+2t^7+t^6+t^5+t^4=-t^4(t-1)^5$

\item \tp{II\, II\,  4 \, 1   } \quad $y^2 =x^3+tx^2+t^2(2t+1)x+t(t^4+2t^3+t^2+2)$\\
$-\Delta =t^9+2t^7+2t^6+t^4=t^4(t-1)^4(t+1)$

\item \tp{II\, II\,  3 \, 2    } \quad $y^2 =x^3+tx^2i+t^3x+t(t^4i+t^3i-t^2i+ti+2i)$ \\
$-\Delta =2t^9+t^8+2t^7+t^6+2t^5+t^4=-t^4(t-1)^3(t+1)^2$

\item \tp{II\, II\,  3 \, 1^2  } \quad $y^2 =x^3+tx^2+t^3x+t(t^4+t^3+2t^2+2)$ \\
$-\Delta =2t^9+2t^7+t^6+t^4=-t^4(t-1)^3\left(t-\frac{1}{i}\right)\left(t+\frac{1}{i}\right)$

\item \tp{II\, II\, 2^2 \, 1    } \quad $y^2 =x^3+tx^2i+t^3x+t(t^4i-t^3+t^2i+ti+2i)$ \\
$-\Delta =2t^9+t^8+t^7+2t^6+2t^5+t^4=-t^4\left(t-\frac{1}{i}\right)^2\left(t+\frac{1}{i}\right)^2(t-1)$

\item \tp{II\, II\,  2 \, 1^3  } \quad $y^2 =x^3+tx^2+t^2(t+1)x+t(2t^4+t^3+2)$\\
$-\Delta =2t^9+2t^8+t^7+t^5=t^4(t-1)^2(2t^3+2t+1)$

\item \tp{II\, II\, 1^5   } \quad $y^2 =x^3+tx^2+t^3x+t(t^4+t^3+t^2+t+2)$ \\
$-\Delta =2t^9+2t^7+2t^6+2t^5+t^4=-t^4(t-1)(t+1)(t^3-t+1)$

\item \tp{II\, III\, 5  } \quad $y^2 =x^3+tx^2i+t^3x+t(t^3i-t^2i-ti+2i)$ \\
$-\Delta =2t^9+2t^8+2t^7+t^6+t^5+t^4=-t^4(t-1)^5$

\item \tp{II\, III\,  4 \, 1   } \quad $y^2 =x^3+t{{\epsilon}}x^2-t^3xi+\frac{t}{{\epsilon}}((2+2i)t^3+(2+i)t^2+t+2+2i)$ \\
$-\Delta =-t^9i+(1-i)t^8+t^7+(3-i)t^6+(1-i)t^5+t^4=-t^4i(t+1)^4\left(t-\frac{1}{i}\right)$

\item \tp{II\, III\,  3 \, 2    } \quad $y^2 =x^3+tx^2+t^3x+t(t^3+2t^2+t+2)$ \\
$-\Delta =2t^9+t^8+2t^7+t^6+2t^5+t^4=-t^4(t-1)^3(t+1)^2$

\item \tp{II\, III\,  3 \, 1^2   } \quad $y^2 =x^3+tx^2+t^2(2t+1)x+t(2t^3+t^2+t+2)$ \\
$-\Delta =t^9+t^8+2t^7+2t^6+2t^5+t^4=-t^4(t-1)^3(2t^2+2t+1)$

\item \tp{II\, III\, 2^2 \, 1    } \quad $y^2 =x^3+tx^2+t^3x+t(2t^3+t^2+t+2)$ \\
$-\Delta =2t^9+t^8+t^7+2t^6+2t^5+t^4=-t^4(t-1)\left(t-\frac{1}{i}\right)^2\left(t+\frac{1}{i}\right)^2$

\item \tp{II\, III\,  2 \, 1^3   } \quad $y^2 =x^3+tx^2+t^2(2t+1)x+t(t^3+2t^2+2t+1)$ \\
$-\Delta =t^9+t^8+t^6+t^5+2t^4=t^4(t-1)^2(t^3+2t+2)$

\item \tp{II\, III\, 1^5  } \quad $y^2 =x^3+tx^2+t^2(t+1)x+t(t^3+2t^2+2)$ \\
$-\Delta =2t^9+t^8+t^7+t^6+t^4=-t^4(t-1)(t+1)(t^3+2t^2+1)$
\end{enumerate}

\subsubsection{Case 1C}\label{sec:1C}
$$\delta=6,\quad r=0\quad \mbox{ Type \tp{II}}$$
\begin{equation}\label{eq.1c}
y^2=x^3+t^2c_0x^2+t^2c_2x+tc_5, \quad t\nmid c_5, t\nmid c_2
\end{equation}
$$-\Delta=t^4c_0^2(t^3c_0c_5-t^4c_2^2)+t^6c_2^3$$

\begin{enumerate}
\itemsep=.5ex
\item \tp{II\, 6  } \quad does not exist,~\ref{algman}

\item \tp{II\,  5 \, 1   } \quad $y^2 =x^3+t^2x^2+t^2(t^3+t+1)x+t(2t^5+t^3+t^2+2)$ \\
$-\Delta =2t^{12}+t^{11}+t^{10}+2t^8+2t^7+t^6=-t^6(t-1)^5(t+1)$

\item \tp{II\,  4 \, 2   } \quad $y^2 =x^3+t^2x^2+t^2(t^2+t+2)x+t(2t^5+t^4+2+2t^2+2t)$ \\
$-\Delta =t^4(t^3(2t^5+t^4+2t^2+2t+2)+2t^4(t^2+t+2)^2+t^6(t^2+t+2)^3=-t^6(t-1)^4(t+1)^2$

\item \tp{II\,  4 \, 1^2  } \quad $y^2 =x^3+t^2x^2+t^2(t^2+t+2)x+t(2t^5+t^3+2t^2+1)$ \\
$-\Delta =2t^{12}+t^{11}+2t^{10}+2t^9+2t^8+t^7+2t^6=-t^6(t-1)^4(t-i)(t+i)$

\item \tp{II\, 3^2  } \quad does not exist,~\ref{algman}

\item \tp{II\, 3\, 2 \, 1   } \quad $y^2 =x^3+t^2x^2+(2t^4+t^3+t^2)x+it^5+t^4+2t^3+2(1+i)t^2+t$ \\
$-\Delta =t^{12}+(2+i)t^{11}+2t^{10}+t^9+(1+2i)t^8+t^7+t^6=-t^6(t-1)^3(-t^3+(1-i)t^2+t+1)$

\item \tp{II\,  3 \, 1^3  } \quad $y^2 =x^3+t^2x^2+t^2(1+t^2+t)x+t(t^5+2t^4+2t^3+2t^2+t+1)$ \\
$-\Delta =t^{12}+2t^{10}+t^9+t^7+t^6=-t^6(t-1)^3(2t^3+t+1)$

\item \tp{II\, 2^3  } \quad $y^2=x^3+2t^2x^2+(t^4+t^3+\phi_1t^2)x+t^2\phi_1^3+2t^5+(2+\phi_1)t^4+(2+\phi_1)t^3+2\phi_1^2t^2+t$\\
$-\Delta =2t^{12}\phi_1^3+2t^{11}+2t^9+2t^7+t^6\phi_1^3=-t^6(t-1)^2(t+1)^2(2t^2\phi_1^3+2t+\phi_1^3)$, where  $\phi_1^6=-1$

\item \tp{II\, 2^2 \, 1^2  } \quad $y^2 =x^3+t^2x^2+(t^4+t^3+2t^2)x+t^6+2t^4+t^3+t^2+t$ \\
$-\Delta =t^{12}+t^{11}+t^9+t^7+2t^6=-t^6(t-1)^2(t+1)^2(2t^2+2t+1)$

\item \tp{II\,  2 \, 1^4  } \quad $y^2 =x^3+t^2x^2+(t^4+t^3+t^2)x+t^6+t^5+t^4+t^3+t^2+t$ \\
$-\Delta =t^{12}+2t^{11}+t^{10}+t^7+t^6=-t^6(t-1)^2(t+1)(2t^3+t+2)$

\item \tp{II\, 1^6  } \quad $y^2 =x^3+t^2x^2+t^2(t^2+t+1)x+t(2t^5+t^3+2t^2+1)$ \\
$-\Delta =t^4(t^3(2t^5+t^3+2t^2+1)+2t^4(t^2+t+1)^2)+t^6(t^2+t+1)^3=-t^6(t-1)(t+1)(t^4+2t^3+t+1)$

\item \tp{II\, II  } \quad $y^2 =x^3+t^2x+t(t^4+t^3+t^2+t+1)$ \\
$-\Delta =2t^6$

\item \tp{II\, II\, II  } \quad $y^2 =x^3+t^2(2t+1)x+t(t^4+t^3+t^2+2t+1)$ \\
$-\Delta =t^9+2t^6=2t^6(2t+1)^3$

\item \tp{II\, II\, III  } \quad $y^2 =x^3+t^2(2t+1)x+t(t^4+t^3+2t^2+t+1)$ \\
$-\Delta =t^9+2t^6=2t^6(2t+1)^3$

\item \tp{II\, III\, III  } \quad $y^2 =x^3+t^2(2t+1)x+t(t^2+t+1)$ \\
$-\Delta =t^9+2t^6=2t^6(2t+1)^3$
\end{enumerate}

\subsubsection{Case 1D}
$$\delta=7,\quad r=0\quad \mbox{ Type \tp{II}}$$
\begin{equation}
y^2=x^3+t^2c_0x^2+t^3c_1x+tc_5, \quad t\nmid c_5, t\nmid c_0
\end{equation}
$$-\Delta =t^4c_0^2(t^3c_0c_5-t^6c_1^2)+t^9c_1^3 $$

\begin{enumerate}
\itemsep=.5ex
\item \tp{II\, 5  } \quad $y^2 =x^3+t^2x^2+t^2(t+1)x+t(2t^5+t^4+t+1)$ \\
$-\Delta =t^4(t^3(2t^5+t^4+t+1)+2t^6(t+1)^2)+t^9(t+1)^3=-t^7(t-1)^5$

\item \tp{II\,  4 \, 1   } \quad $y^2 =x^3+t^2x^2+t^2(t+1)x+t(2t^5+t^4+2t^3+2)$ \\
$-\Delta =2t^{12}+t^{10}+t^9+2t^7=-t^7(t+1)(t-1)^4$

\item \tp{II\,  3 \, 2   } \quad $y^2 =x^3+t^2x^2+t^2(t+1)x+t(2t^5+2t+1)$ \\
$-\Delta =2t^{12}+t^{11}+2t^{10}+t^9+2t^8+t^7=-t^7(t-1)^3(t+1)^2$

\item \tp{II\,  3 \, 1^2  } \quad $y^2 =x^3+2t^2x^2+(t^3+t^2)x+t^6+t^4+t^2+t$ \\
$-\Delta =2t^{12}+t^{11}+t^{10}+t^9+2t^8+2t^7 =-t^7(t-1)^3(t^2+2t+2)$

\item \tp{II\,  2 \, 1^3   } \quad $y^2 =x^3+2t^2x^2+t^2(t+1)x+t(2t^5+t^2+1)$ \\
$-\Delta =t^4(2t^3(2t^5+t^2+1)+2t^6(t+1)^2)+t^9(t+1)^3=-t^7(t-1)^2(t+1)(2t^2+t+1)$

\item \tp{II\, 2^2 \, 1    } \quad $y^2 =x^3+t^2x^2+t^2(t+1)x+t(2t^5+2t^3+t^2+2t+1)$ \\
$-\Delta =2t^{12}+t^{11}+t^{10}+2t^9+2t^8+t^7=-t^7(t-i)^2(t+i)^2(t-1)$

\item \tp{II\, 1^5   } \quad $y^2 =x^3+t^2x^2+t^2(t+1)x+t(t^5+2t^4+2t^3+t^2+t+1)$ \\
$-\Delta =t^4(t^3(t^5+2t^4+2t^3+t^2+t+1)+2t^6(t+1)^2)+t^9(t+1)^3=-t^7(t-1)(t+1)(2t^3+t+1)$

\item \tp{II\, IV  } \quad  does not exist,~\ref{b2}

\item \tp{II\, II\, 1  } \quad does not exist,~\ref{b2}

\item \tp{II\, II\, 2  } \quad does not exist,~\ref{b2}

\item \tp{II\, II\, 1^2  } \quad does not exist,~\ref{b2}

\item \tp{II\, III\, 2   } \quad does not exist,~\ref{b2}

\item \tp{II\, III\, 1^2  } \quad does not exist,~\ref{b2}
\end{enumerate}

\subsubsection{Case 1E}
$$\delta=9,\quad r=0\quad \mbox{ Type \tp{II}}$$
\begin{equation}
y^2=x^3+t^3c_1x+tc_5, \quad t\nmid c_5, t\nmid c_1
\end{equation}
$$-\Delta=t^9c_1^3$$

\begin{enumerate}
\itemsep=.5ex
\item \tp{II\, 3  } \quad does not exist,~\ref{b2}

\item \tp{II\,  2 \, 1   } \quad does not exist,~\ref{b2}

\item \tp{II\, 1^3  } \quad does not exist,~\ref{b2}

\item \tp{II\, II   } \quad $y^2 =x^3+t^3(2t+1)x+t(t^3+t+1)$ \\
$-\Delta =t^9(2t+1)^3$

\item \tp{II\, III  } \quad $y^2 =x^3+t^3(2t+1)x+t(t^2+t+1)$ \\
$-\Delta =t^9(2t+1)^3$
\end{enumerate}

\subsubsection{Case 1F}
$$\delta=12,\quad r=0\quad \mbox{ Type \tp{II }}$$
\begin{equation}
y^2=x^3+t^4c_0x+tc_5, \quad t\nmid c_5, t\nmid c_0
\end{equation}
$$-\Delta=t^{12}c_0^3$$

\subsubsection{Case 2}
$$\delta=3,\quad r=1\quad \mbox{ Type \tp{III }}$$
\begin{equation}y^2=x^3+tc_1x^2+tc_3x+t^2c_4, \quad t\nmid c_3
\end{equation}
$$-\Delta=t^2c_1^2(t^3c_1c_4-t^2c_3^2)+t^3c_3^3$$

\begin{enumerate}
\itemsep=.5ex
\item \tp{III\, III\, III\, III  } \quad 
$y^2 = x^3-t(t-1)(t-(\alpha+\alpha^3))x+ (\alpha^5-\alpha^3)(t^2(t-1)^2)$\\
$-\Delta =t^9-t^6(1+(\alpha^3-\alpha)) + (\alpha^3-\alpha)t^3 = t^3(t-1)^3(t-(\alpha+\alpha^3))^3$,\\
where $\alpha^4+\alpha^2 = 1$.

\medskip

\item \tp{III\, III\, III\, II  } \quad 
$y^2 = x^3-t(t^2-1)x+ t^2(t-1)^2$\\
$-\Delta =t^9+t^3 = t^3(t^2+1)^3$

\item \tp{III\, III\, II\, II  } \quad 
$y^2 = x^3-t(t^2+1)x+ t^2(t^2+1)^2$\\
$-\Delta =t^9+t^3 = t^3(t^2+1)^3$

\item \tp{III\, II\, II\, II  } \quad 
$y^2 = x^3-t(t^2+1)x+ t^2(t^2+1)(t-1)$\\
$-\Delta =t^9+t^3 = t^3(t^2+1)^3$

\item \tp{III\, II\, 6  } \quad does not exist,~\ref{algman}

\item \tp{III\, II\,  5 \, 1   } \quad
$y^2 =x^3+tx^2i+t(t^2+2)x+t^2(t^3i-t^2i-ti+i)$ \\
$-\Delta =2t^9+t^8+t^7+2t^5+2t^4+t^3$

\item \tp{III\, II\,  4 \, 2   } \quad $y^2 =x^3+tx^2i+t(t^2+1)x+t^2(t^3+t^2i+(1-i)t+2+2i)$ \\
$-\Delta =2t^9-(1-i)t^8-t^7+(1-i)t^6+(1-i)t^5+t^4+t^3i$

\item \tp{III\, II\,  4 \, 1^2  } \quad
$y^2 =x^3+tx^2+t(2t^2+1)x+t^2(2t^3+t^2+t+2)$ \\
$-\Delta =t^9+2t^8+2t^7+t^5+t^4+2t^3 =t^3(t-1)^4(t^2+2)$

\item \tp{III\, II\, 3^2  } \quad does not exist,~\ref{algman}

\item \tp{III\, II\, 3\, 2 \, 1  } \quad $y^2 =x^3+t{{\epsilon}}x^2+t(2t^2+2)x+\frac{t^2}{{\epsilon}}((1-i)t^3+t^2+it+i)$ \\
$-\Delta = t^9+(1+i)t^8+(2+i)t^7+t^6+(2-i)t^5+(1-i)t^4+t^3 =t^3(t-1)^3(t-i)^2(t+1))$

\item \tp{III\, II\,  3 \, 1^3  } \quad
$y^2 =x^3+tx^2+t(t^2+1)x+t^2(t^3+t^2)$ \\
$-\Delta =2t^9+2t^7+2t^6+t^4+2t^3=-t^3(t-1)^3(t^3+t+2)$

\item \tp{III\, II\, 2^3  } \quad
$y^2 =x^3+tx^2+t((1-i)t^2+1+i)x+t^2(-it^3+2t^2+(1+i)t+2)$ \\
$-\Delta =(-1-i)t^9-it^8+t^7-it^6+t^5+2it^4+(-1+i)t^3=t^3(t-1)^2(t+1)^2(-1+i-it+(-1-i)t^2)$

\item \tp{III\, II\, 2^2 \, 1^2  } \quad
$y^2 =x^3+tx^2+t(t^2i+1)x+t^2(t^3+(1-i)t^2-(1+i)t+1-i)$ \\
$-\Delta =t^9i+t^8+(2+i)t^7+t^6-(1+2i)t^5+t^4-t^3=t^3(t-1)^2(t+1)^2(2+t+it^2)$

\item \tp{III\, II\,  2 \, 1^4  } \quad
$y^2 = x^3+tx^2+t(t^2+1)x+t^2(t^3+2+2t)$ \\
$-\Delta =2t^9+t^5+t^4+2t^3=t^3(t-1)^2(2+2t+t^3+2t^4)$

\item \tp{III\, II\, 1^6  } \quad
$y^2 =x^3+tx^2+t(2t^2+2)x+t^2(2t^3+t^2+2t+1)$ \\
$-\Delta =t^9+2t^8+2t^7+2t^5+t^4+t^3=t^3(t-1)(t+1)(t-i)(t+i)\left(t-\frac{1}{1-i}\right)\left(t-\frac{1}{1+i}\right)$

\item \tp{III\, III\, 6  } \quad does not exist,~\ref{algman}

\item \tp{III\, III\,  5 \, 1   } \quad
$y^2 =x^3+t{{\epsilon}}x^2+t^3xi+tx+\frac{t^4}{{\epsilon}}+\frac{t^3}{(2+2i){\epsilon}}+\frac{2t^3i}{2+2i){\epsilon}}+\frac{2t^2}{{\epsilon}}$ \\
$-\Delta=t^9i-(2+2i)t^8+(1+1i)t^7+\left(\frac{1}{1+i}+2+i\right)t^6+(2+2i)t^5+(2+2i)t^4+2t^3$, where ${\epsilon}^2+1+i=0$

\item \tp{III\, III\,  4 \, 2   } \quad does not exist,~\ref{algman}

\item \tp{III\, III\,  4 \, 1^2  } \quad $y^2=x^3+t{{\epsilon}}x^2+t(t^2+1+i)x+\frac{t^2}{{\epsilon}}(-t^2i-\left(\frac{1}{2}-\frac{1i}{2}\right)t+1)$\\
$-\Delta =2t^9-(1+i)t^8+(1-i)t^7-(1+i)t^6-(2+2i)t^5+(2+i)t^4+(2+i)t^3$, where ${\epsilon}^2+1+i=0$

\item \tp{III\, III\, 3^2  } \quad does not exist,~\ref{algman}

\item \tp{III\, III\, 3\, 2 \, 1  } \quad $y^2=x^3+t{\epsilon}x^2+t((1-i)t^2+2)x+\frac{t^2}{{\epsilon}}(-t^2+\left(\frac{1}{5}-\frac{2i}{5}\right)t+1)$\\
$-\Delta =(2+2i)t^9+(2+2i)t^8+2it^7+t^6+2it^5-(1-i)t^4+t^3$

\item \tp{III\, III\,  3 \, 1^3  } \quad $y^2 =x^3+tx^2+t(2t^2+1)x+t^2(1+t^2+t)$ \\
$-\Delta =t^9+t^8+2t^7+2t^5+t^4+2t^3=-t^3(t-1)^3(2t^3+2t^2+t+2)$

\item \tp{III\, III\, 2^3  } \quad $y^2 =x^3+tx^2+t(2t^2+t+2)x+t^2(1+t^2)$ \\
$-\Delta =t^9+t^8+2t^6+t^4+t^3 =t^3(t-1)^2(t-i)^2(t+i)^2$

\item \tp{III\, III\, 2^2 \, 1^2   } \quad $y^2 =x^3+tx^2+t(t^2+2)x+t^2(t^2+2t+2)$ \\
$-\Delta =2t^9+t^8+2t^7+2t^6+t^5+t^4+t^3 =-t^3\left(\frac{t-1}{i}\right)^2\left(\frac{t+1}{i}\right)^2
\left(t-\frac{1}{1+i}\right)\left(t-\frac{1}{1-i}\right)$

\item \tp{III\, III\,  2 \, 1^4  } \quad $y^2 =x^3+tx^2+t(2t^2+1)x+t^2(t^2i+ti+i)$ \\
$-\Delta =t^9+2t^8+2t^7+t^6+2t^5+2t^4+2t^3 =t^3(t-1)^2(-1+t^3+t^4)$

\item \tp{III\, III\, 1^6  } \quad
$y^2 =x^3+tx^2+t(t^2+t+1)x+t^2(1+t^2)$ \\
$-\Delta =2t^9+t^8+2t^7+t^5+t^4+2t^3 =-t^3(t-1)(t^5+t^3+t^2+2)$

\item \tp{III\, 9  } \quad does not exist,~\ref{sigma-r}

\item \tp{III\,  8 \, 1   } \quad does not exist,~\ref{add3}

\item \tp{III\,  7 \, 2   } \quad does not exist,~\ref{add3}

\item \tp{III\,  7 \, 1^2   } \quad
$y^2 =x^3+t(t+1)x^2+t(2t^2+2)x+t^4$ \\
$-\Delta =2t^5+t^4+t^3  =t^3(1+t-t^2)$

\item \tp{III\,  6 \, 3   } \quad
$y^2 =x^3+t^2x^2+t(t^2+1)x+t^2(2+t^2+2t)$ \\
$-\Delta =t^6+2t^3  =t^3(t-1)^3$

\item \tp{III\, 6\, 2 \, 1   } \quad does not exist,~\ref{add3}

\item \tp{III\,  6 \, 1^3  } \quad
$y^2 =x^3+t(t+1)x^2+t(t^2+1)x+t^2(t^2+t)$ \\
$-\Delta =2t^6+2t^5+t^4+2t^3 =t^3(2+t-t^2-t^3)$

\item \tp{III\,  5 \, 4   } \quad does not exist,~\ref{add3}

\item \tp{III\, 5\, 3 \, 1   } \quad does not exist,~\ref{add3}

\item \tp{III\,  5 \, 2^2  } \quad does not exist,~\ref{lattice}

\item \tp{III\, 5\, 2 \, 1^2  } \quad
$y^2 =x^3+tx^2+t(-it^3+2)x+t((1+i)t^5+2t^4+(1-i)t^3+t^2i+t)$\\
$-\Delta =2it^{12}+2t^{10}+(-1-i)t^9+t^8+2t^7+2it^6+2t^5+t^4+t^3$

\item \tp{III\,  5 \, 1^4  } \quad
$y^2 =x^3+t(t+1)x^2+t(2t^2+i)x+t^2(t^2+1+i)$ \\
$-\Delta =(-1-i)t^7+(2+i)t^6-it^5+2t^4+t^3i =-t^3(t-1)(i-(1-i)t-t^2+(1+i)t^3)$

\item \tp{III\, 4^2 \, 1   } \quad does not exist~\ref{algman}

\item \tp{III\, 4\, 3 \, 2   } \quad does not exist,~\ref{add3}

\item \tp{III\, 4\, 3 \, 1^2  } \quad
$y^2 =x^3+tx^2+t((-1-i)t^3+1+2i)x+t(2t^3+(1-i)t^2)$\\
$-\Delta =(1-i)*t^{12}+2it^{10}+(2+i)t^6-2it^4-(1+i)t^3$

\item \tp{III\, 4\, 2^2 \, 1   } \quad
$y^2 =x^3+t(t+1)x^2+t((-1-i)t^2+2+2i)x+t^2(-it^2+t+i)$ \\
$-\Delta =-it^8+2t^6+2it^4+(1-i)t^3  =-t^3(t-1)^2\left(t+\frac{1}{i}\right)^2(it+(1-i))$

\item \tp{III\, 4\, 2 \, 1^3  } \quad
$y^2 =x^3+t(t+1)x^2+t(t^2+1)x+t^2(2+t^2+t)$ \\
$-\Delta =t^8+2t^6+t^4+2t^3 =t^3(t-1)^2(2-t-t^2+t^3)$

\item \tp{III\,  4 \, 1^5  } \quad
$y^2 =x^3+t(t+1)x^2+t(t^2+1)x+t^2(1+t^2+t)$ \\
$-\Delta =2t^8+2t^6+t^5+t^4+2t^3 =t^3(2+t+t^2-t^3-t^5)$

\item \tp{III\, 3^3  } \quad
$y^2 =x^3+t^2x^2+t(t^2+2)x+t^2(1+t^2)$ \\
$-\Delta =2t^9+t^6+t^3 =-t^3\left(t-\frac{1}{(1-i)}\right)^3\left(t-\frac{1}{(1+i)}\right)^3$

\item \tp{III\, 3^2\, 2 \, 1   } \quad does not exist,~\ref{add3}

\item \tp{III\, 3^2 \, 1^3  } \quad
$y^2 = x^3+(t^2+ti)x^2+(\phi_0t+t^3i)x+2t^4+t^3\phi_0 $\\
$- =(2i+2\phi_0)t^9+(-i\phi_0+1)t^8+(2i+2\phi_0)t^7+(2i+2\phi_0-i\phi_0)t^6+(\phi_0i+2)t^5+(i+\phi_0)t^4+(2\phi_0+2i+\phi_0i)t^3 =-t^3(t-1)^3((\phi_0+i)t^3+(-1+i\phi_0)t^2+(i+\phi_0)t+(-\phi_0-i+i\phi_0))$, where $\phi_0^2+\phi_0+i=0$

\item \tp{III\,  3 \, 2^3  } \quad
$y^2 =x^3+t(t+1)x^2+2tx+t^2(2+2t)$ \\
$-\Delta =t^9+t^8+2t^6+t^4+t^3  =t^3(t-1)^2(t-i)^2(t+i)^2$

\item \tp{III\, 3\, 2^2 \, 1^2  } \quad
$y^2 = x^3+tx^2+t((1+2i)t^3+1+i)x+t(t^5i+t^4i+t^3+(1-i)t^2-it)$ \\
$-\Delta = ((-1-i)t^{12}+it^{10}+2it^9+2it^8+(2+i)t^6+it^5+2it^4+(-1+i)t^3= t^3(t-1)^3(t^6+(1+i)t^4-i)t^3-(1+2i)t^2+(1+i)t+i)$

\item \tp{III\, 3\, 2 \, 1^4  } \quad
$y^2 =x^3+t(t+1)x^2+t((1+i)t^2+1)x+t^2(-it^2+(1+i)t+1+i)$ \\
$-\Delta =(1+i)t^9-(2+3i)t^8+(1+2i)t^7+(-1+i)t^6+(4-i)t^5+t^4+2t^3=t^3(t-1)^2(2-t-it^2-it^3+(1+i)t^4)$

\item \tp{III\,  3 \, 1^6  } \quad
$y^2 =x^3+t^2x^2+tx+t^2(1+t)$ \\
$-\Delta =t^9+t^8+-t^6+t^3 =t^3(t^6+t^5-t^3+1)$

\item \tp{III\, 2^4 \, 1   } \quad
$y^2 =x^3+t{\epsilon}x^2+t(2t^3+2+i)x+\frac{t^2}{{\epsilon}}((-1+i)t^4-(1+i)t^3-(1+i)t^2+ti+i)$\\
$-\Delta = t^{12}+(2-i)t^{10}+t^9+4it^8+(2+4i)t^7+(2+4i)t^6+(2+4i)t^5-(2+4i)t^4+(4+4i)t^3=t^3(t-(-1-i))^2(t-1)^2(t+1)^2(t+i)^2(t-(-1+i))$, where ${\epsilon}^2=2-i$

\item \tp{III\, 2^3 \, 1^3  } \quad
$y^2 =x^3+t^2x^2+tx+t^2(2+2t^2+2t)$ \\
$-\Delta =t^{10}+t^9+t^8+t^6+2t^3 =t^3(t-1)^2(t+1)^2(2+t^2+t^3)$

\item \tp{III\, 2^2 \, 1^5  } \quad
$y^2 =x^3+t^2x^2+t(t^2+2)x+t^2(2+2t^2+2t)$ \\
$-\Delta =2t^{10}+2t^8+t^6+t^3 =t^3(t-1)^2(1-t-t^3+t^4-t^5)$

\item \tp{III\,  2 \, 1^7  } \quad
$y^2 =x^3+(t+t^2)x^2+2tx+t^4+2t^3$ \\
$-\Delta =2t^{10}+t^9+2t^7+2t^6+2t^5+t^4+t^3=t^3(t-i)(t+i)(1+t+t^2+t^3+t^4-t^5)$

\item \tp{III\, 1^9  } \quad
$y^2=x^3+t(t+1)x^2+2tx+t^2(t^3+2)$ \\
$-\Delta =2t^{11}+t^6+t^4+t^3 =t^3(t-i)(t+i)(1+t-t^2+t^4-t^6)$

\end{enumerate}

\subsubsection{Case 3A}
$$\delta=5,\quad r=2\quad \mbox{ Type \tp{IV }}$$
\begin{equation}
y^2=x^3+tc_1x^2+t^2c_2x+t^2c_4, \quad t\nmid c_4, t\nmid c_1
\end{equation}
$$-\Delta=t^2c_1^2(t^3c_1c_4-t^4c_2^2)+t^6c_2^3$$

\begin{enumerate}
\item \tp{IV\, 7  } \quad does not exist,~\ref{lattice}

\item \tp{IV\,  6 \, 1   } \quad does not exist,~\ref{add3}

\item \tp{IV\,  5 \, 2   } \quad does not exist,~\ref{algman}

\item \tp{IV\,  5 \, 1^2   } \quad
$y^2 =x^3+t(t+1)x^2+2t^2x+t^2$ \\
$-\Delta =2t^7+2t^5 =t^5(2t^2+2)^2$

\item \tp{IV\,  4 \, 3   } \quad does not exist,~\ref{add3}

\item \tp{IV\, 4\, 2 \, 1   } \quad
$y^2 =x^3+tx^2i+t^2(x^2+x+2xi+2)$ \\
$-\Delta =(1+i)t^8+t^7+t^6-it^5$

\item \tp{IV\,  4 \, 1^3  } \quad
$y^2 =x^3+t(t+1)x^2+2t^2xi-t^2(2+i)$ \\
$-\Delta =-(2+2i)t^8+t^7-(1+2i)t^6-(1+2i)t^5$

\item \tp{IV\, 3^2 \, 1   } \quad
$y^2 =x^3+tx^2i+t^2(x^2+2i)-t^3i$ \\
$-\Delta =t^9i+t^8i+t^6+t^5$

\item \tp{IV\,  3 \, 2^2   } \quad
$y^2 =-\frac{1}{2}t^3i+(2+2i)(-2xi+x-x^2i+2+2i+x^2)t^2 +(2+2i)(x^2i+x^2)t+(2+2i)(x^3-x^3i)$ \\
$-\Delta =-it^9 -(1-i)t^8+t^7i-it^6+t^5$

\item \tp{IV\, 3\, 2 \, 1^2  } \quad
$y^2 =x^3+tx^2+t^2(x^2+2x+1)+2t^3$ \\
$-\Delta =t^9+2t^7+t^6+2t^5$

\item \tp{IV\,  3 \, 1^4  } \quad
$y^2 =x^3+t(t+1)x^2+2t^2x+t^2(t+1)$ \\
$-\Delta =2t^9+2t^7+2t^6+2t^5$

\item \tp{IV\, 2^3 \, 1   } \quad
$y^2 =x^3+t(t+1)x^2+2t^2x+t^2((2+i)t^2+(2+2i)t+2+2i)$ \\
$-\Delta=(1+2i)t^{10}+(1+i)t^9+(2+i)t^8+(2+i)t^7+(2-i)t^6+(1-i)t^5$

\item \tp{IV\, 2^2 \, 1^3  } \quad
$y^2 =x^3+t(t+1)x^2+2t^2x+t^2(t^3\gamma_3+2t^2+t)$ \\
$-\Delta =t^{10}+2t^9+t^8+2t^6$

\item \tp{IV\,  2 \, 1^5  } \quad
$y^2 =x^3+t(t+1)x^2+t^2x+t^2(t^2+2t+1)$ \\
$-\Delta =2t^{10}+t^9+t^7+2t^5$

\item \tp{IV\, 1^7  } \quad
$y^2 =x^3+t(t+1)x^2+t^2x+t^2(2t^3+t^2+1)$ \\
$-\Delta =t^{11}+2t^{10}+t^8+t^7+2t^6+2t^5$

\item \tp{IV\, IV\, 2  } \quad
$y^2 =x^3+tx^2+2t^2x+t^2(2t^2+2t+2)$ \\
$-\Delta =t^7+t^6+t^5 =t^5(t-1)^2$

\item \tp{IV\, IV\, 1^2  } \quad
$y^2=x^3+tx^2+2t^2x+t^2(2+t^2) $\\
$-\Delta =2t^7+t^5 =2t^5(t-1)(t+1)$

\item \tp{IV\, II\, 3  } \quad
$y^2 =x^3+tx^2+2t^2x+t^2(t^3+2)$ \\
$-\Delta =2t^8+t^5=t^5(t-1)^2$

\item \tp{IV\, II\,  2 \, 1   } \quad
$y^2 =x^3+tx^2+2t^3\phi_1x+t^2(2t^3+t^2+t+2)$ \\
$-\Delta =t^8+2t^7+2t^6+t^5$

\item \tp{IV\, II\, 1^3  } \quad
$y^2 =x^3+tx^2+t^2(2t^3+t+2)$ \\
$-\Delta =t^8+2t^6+t^5$

\item \tp{IV\, III\, 4  } \quad
$y^2 =x^3+tx^2i+t^3x+t^2(ti+2i)$ \\
$-\Delta =t^9+2t^8+2t^6+t^5 =t^5(t-1)^4$

\item \tp{IV\, III\,  3 \, 1   } \quad
$y^2 =x^3+tx^2i+2t^3x+t^2(2ti+2i)$ \\
$-\Delta =2t^9+2t^8+t^6+t^5 =-t^5(t-1)^3(t+1)$

\item \tp{IV\, III\, 2^2  } \quad
$y^2 =x^3+2tx^2+2t^2(1+i)t+1)x+t^2((2+2i)t^2+2t+1+2i)$ \\
$-\Delta =-(1-i)t^9+2it^8-(2+2i)t^7+2t^6-(2+i)t^5$

\item \tp{IV\, III\,  2 \, 1^2  } \quad
$y^2=x^3+tx^2+2t^2(2t+1)x+t^2(2t^2+2t+2)$ \\
$-\Delta =t^9+t^8+2t^7+t^6+t^5=t^5(t-1)^2(t-i)(t+i)$

\item \tp{IV\, III\, 1^4  } \quad
$y^2=x^3+tx^2+2t^2(t+2)x+t^2(t^2+t+2)$ \\
$-\Delta =2t^9+t^8+t^6+t^5=t^5(2t^4+t^3+t+1)$

\item \tp{IV\, II\, 4  } \quad
$y^2 =x^3+tx^2+2t^2(2t+1)x+t^2(2t^3+t^2+t+2)$ \\
$-\Delta =t^9+2t^8+2t^6+t^5=t^5(t-1)^4$

\item \tp{IV\, II\,  3 \, 1   } \quad $y^2 =x^3+tx^2+2t^2(t+1)x+t^2(2t^3+2t^2+2t+2)$ \\
$-\Delta =2t^9+2t^8+t^6+t^5=t^5(t-1)^3(t+1)$

\item \tp{IV\, II\, 2^2  } \quad $y^2 =x^3+tx^2+2t^2(2t+1)x+t^2(t^3+2)$ \\
$-\Delta =t^9+t^7+t^5=t^5(t-1)^2(t+1)^2$

\item \tp{IV\, II\,  2 \, 1^2  } \quad $y^2 =x^3+tx^2i+t^3x+t^2(t^3i+t^2i+2ti+2i)$ \\
$-\Delta =t^9+t^8+2t^7+t^6+t^5=t^5(t-1)^2(t-i)(t+i)$

\item \tp{IV\, II\, 1^4  } \quad $y^2 =x^3+tx^2+2t^2(t+1)x+t^2(2t^3+2t^2+t+2)$ \\
$-\Delta =2t^9+2t^8+2t^6+t^5=t^5(2t^4+2t^3+2t+1)$

\end{enumerate}

\subsubsection{Case 3B}\label{sec:3B}
$$\delta=6,\quad r=2\quad \mbox{ Type \tp{IV }}$$
\begin{equation}
y^2=x^3+t^2c_0x^2+t^2c_2x+t^2c_4, \quad t\nmid c_4, t\nmid c_2
\end{equation}
$$-\Delta=t^4c_0^2(t^4c_0c_4-t^4c_2^2)+t^6c_2^3$$

\begin{enumerate}
\itemsep=.5ex
\item \tp{IV\, 6  } \quad $y^2 =x^3+t^2x^2+(t^4+t^3+2t^2)x+2t^6+2t^5+2t^3+t^2$ \\
$-\Delta =2t^{12}+2t^9+2t^6=2t^6(t^6+t^3+1)$

\item \tp{IV\,  5 \, 1   } \quad does not exist,~\ref{algman}

\item \tp{IV\,  4 \, 2   } \quad does not exist,~\ref{algman}

\item \tp{IV\,  4 \, 1^2  } \quad $y^2 =x^3+t^2x^2+(t^4+t^3+t^2)x+2t^6+t^5+t^3$ \\
$-\Delta =2t^{12}+2t^{11}+2t^8+t^6=t^6((t-1)^4(t-(-1-i))(t-(-1+i))$

\item \tp{IV\, 3^2  } \quad does not exist,~\ref{algman}

\item \tp{IV\, 3\, 2 \, 1   } \quad does not exist,~\ref{algman}

\item \tp{IV\,  3 \, 1^3  } \quad $y^2 =x^3+2t^2x^2+(t^3+t^2)x+2t^6+t^5+2t^4+t^3+t^2$ \\
$-\Delta =2t^{12}+2t^{11}+t^9+t^8+t^6=-t^6(t^3+2)(2t^3+t^2+1)$

\item \tp{IV\, 2^3  } \quad does not exist,~\ref{algman}

\item \tp{IV\, 2^2 \, 1^2  } \quad $y^2 =x^3+t^2x^2+(t^3+2t^2)x+2t^6+2t^4+2t^2$ \\
$-\Delta =2t^{12}+t^{10}+t^8+2t^6=t^6(1+t^2)(2t^4+2t^2+2)$

\item \tp{IV\,  2 \, 1^4  } \quad $y^2 =x^3+t^2x^2+t^2(t^2+t+1)x+t^2(2t^4+t^2+2t+1)$ \\
$-\Delta =2t^6+t^5+t^4+t^3+1=-t^6(t-1)^2(t-i)(t-(1-i))(t-(1+i))$

\item \tp{IV\, 1^6  } \quad $y^2 =x^3+t^2x^2i+t^2(t^2+t+1)x+t^2(3it^4-2it^3-2it^2+i)$ \\
$-\Delta =2t^{12}+t^{10}+2t^8+t^6=-t^6(t-1)(t+1)(t^4+1)$

\item \tp{IV\, II  } \quad $y^2 =x^3+t^2x+t^2(t^3+1+t)$ \\
$-\Delta =2t^6$

\item \tp{IV\, II\, II  } \quad $y^2 =x^3+t^2(2t+1)x+t^2(t^3+t+1)$ \\
$-\Delta =t^9-t^6=2t^6(2t+1)^3$

\item \tp{IV\, III\, II  } \quad $y^2 =x^3+t^2(2t+1)x+t^2(2t+1)$ \\
$-\Delta =t^9-t^6=2t^6(2t+1)^3$

\item \tp{IV\, III\, III  } \quad $y^2 =x^3+t^2(2t+1)x$ \\
$-\Delta =t^9-t^6=2t^6(2t+1)^3$

\item \tp{IV\, IV  } \quad $y^2 =x^3+2t^2x+t^2(t^2+t+1)$ \\
$-\Delta =t^6$

\end{enumerate}

\subsubsection{Case 3C}
$$\delta=8,\quad r=2\quad \mbox{ Type \tp{IV }}$$
\begin{equation}
y^2=x^3+t^2c_0x^2+t^3c_1x+t^2c_4, \quad t\nmid c_4, t\nmid c_0
\end{equation}
$$-\Delta=t^4c_0^2(t^4c_0c_4-t^6c_1^2)+t^9c_1^3$$

\begin{enumerate}
\itemsep=.5ex
\item \tp{IV\, 4  } \quad $y^2=x^3+2t^2x^2+t^2(t^4+2t^3+2t+1)$\\
$\Delta =t^7(t^5+2t^4+2t^2+t)=t^{12}+2t^{11}+2t^9+t^8$

\item \tp{IV\,  3 \, 1   } \quad $y^2 =x^3+t^2x^2+t^3x+t^2(2t^4+2t^3+t^2+1)$\\
$\Delta =t^{12}+t^{11}+2t^8+2t^9$

\item \tp{IV\, 2^2  } \quad $y^2 =x^3+t^2x^2+t^3(t+2)x+t^2(2t^4+t^3+2+1)$\\
$\Delta =t^{12}+t^8+t^{10}$

\item \tp{IV\,  2 \, 1^2  } \quad $y^2 =x^3+t^2x^2+t^3(t+1)x+t^2(2t^4+t^3+2t^2+t+2)$\\
$\Delta =2t^7(2t^5+t)=t^{12}-t^8$

\item \tp{IV\, 1^4  } \quad $y^2 =x^3+t^2x^2+t^2(2t^4+1)$ \\
$-\Delta =t^8(2t^4+2t^3+t^2+2t+2)$

\item \tp{IV\, II  } \quad does not exist,~\ref{b2}

\item \tp{IV\, II\, 1  } \quad does not exist,~\ref{b2}

\item \tp{IV\, III\, 1  } \quad does not exist,~\ref{b2}
\end{enumerate}

\subsubsection{Case 3D}
$$\delta=9,\quad r=2\quad \mbox{ Type \tp{IV }}$$
\begin{equation}
y^2=x^3+t^3c_1x+t^2c_4, \quad t\nmid c_4, t\nmid c_1
\end{equation}
$$-\Delta=t^9c_1^3$$

\begin{enumerate}
\item \tp{IV\, 3  } \quad does not exist,~\ref{b2}

\item \tp{IV\,  2 \, 1   } \quad does not exist,~\ref{b2}

\item \tp{IV\, 1^3   } \quad does not exist,~\ref{b2}

\item \tp{IV\, II  } \quad $y^2 =x^3+t^3(2t+1)x+t^2$ \\
$-\Delta =t^9(2t+1)^3$

\item \tp{IV\, III  } \quad
$y^2 =x^3+t^3(2t+1)x+t^2(t^2+t+1)$ \\
$-\Delta =t^9(2t+1)^3$

\end{enumerate}

\subsubsection{Case 3E}
$$\delta=12,\quad r=2\quad \mbox{ Type \tp{IV }}$$
\begin{equation}
y^2=x^3+t^4c_0x+t^2c_4, \quad t\nmid c_4, t\nmid c_0
\end{equation}
$$-\Delta=t^{12}c_0^3$$

\subsubsection{Case 4}

In \cite{Lang-II} there are two different sub-cases to this case, but both have Weierstrass form and discriminant as listed below, and thus we may treat them together.

$$\delta=6,\quad r=4\quad \mbox{Type \tp{I_0^*}} $$
\begin{equation}
y^2=x^3+tc_1x^2+t^2c_2x+t^3c_3
\end{equation}
$$-\Delta=t^2c_1^2(t^4c_1c_3-t^8c_0^2)+t^{12}c_0^2$$

\begin{enumerate}
\itemsep=.5ex
\item \tp{ {I_0^*\, 6}  } \quad  does not exist,~\ref{sigma-r}

\item \tp{ {I_0^*\,  5 \, 1 }  } \quad does not exist,~\ref{lattice}

\item \tp{ {I_0^*\,  4 \, 2 }  } \quad does not exist,~\ref{lattice}

\item \tp{ {I_0^*\,  4 \, 1^2}  } \quad $y^2 =x^3+t^2x^2+2t^2x+t^3(t^3+t^2)$ \\
$-\Delta =t^4(t^5(t^3+t^2)+2t^4)+2t^6=t^6(t^6+t^5+2t^2+2)$

\item \tp{ {I_0^*\, 3^2}  } \quad  does not exist,~\ref{lattice}

\item \tp{ {I_0^*\, 3\, 2 \, 1 }  } \quad does not exist,~\ref{lattice}

\item \tp{ {I_0^*\,  3 \, 1^3}   } \quad $y^2 =x^3+tx^2+t^4x+t^3(t^3+t+1)$ \\
$-\Delta =t^2(t^4(t^3+t+1)+2t^8)+t^{12}=t^6(t^3+2t+2)(t^3+2)$

\item \tp{ {I_0^*\, 2^3}   } \quad $y^2 =\left(\frac{1}{{\epsilon}}+\frac{i}{{\epsilon}}\right)t^6-\frac{it^5}{{\epsilon}}\left(2x+\frac{1-i}{{\epsilon}}\right)t^4-\frac{t^3}{{\epsilon}}+tx^2{\epsilon}+x^3$ \\
$-\Delta =t^{12}+it^{10}+(1-i)t^9+2t^8-(1+i)t^7-2it^6=t^6(t+1+i)^2(t-1)^2(t-i)^2$

\item \tp{ {I_0^*\, 2^2 \, 1^2}  } \quad $y^2 =x^3+t^2x^2i+t^2(t^2+t+2)x+t^3(t^2i-ti+i)$ \\
$-\Delta =2t^{12}+t^{10}+t^8+2t^6=t^6(2t^6+t^4+t^2+2)$

\item \tp{ {I_0^*\,  2 \, 1^4}   } \quad $y^2 =x^3+tx^2+t^2(t^2+t+1)x+t^3(2t^2+2t+2)$ \\
$-\Delta =-t^{12}+t^{10}+t^9+t^8+t^6=t^6(t-1)^2(-t^4+t^3+t^2-t+1)$

\item \tp{ {I_0^*\, 1^6}  } \quad $y^2 =x^3+tx^2+t^2(2t^2+2t+1)x+t^3(t^2+t+2)$ \\
$-\Delta =t^{12}+t^{10}+t^8+t^6=t^6(t-i)(t+i)(t^4-1)$

\item \tp{ {I_0^*\, I_0^*}  } \quad $y^2 =x^3+2t2x+t^3$ \\
$-\Delta =t^6$

\item \tp{ {I_0^*\, IV}  } \quad $y^2 =x^3+2t^2x+t^3(t+1)$ \\
$-\Delta =t^6$

\item \tp{ {I_0^*\, II}   } \quad $y^2 =x^3+2t^2x+t^3(t^2+t+1)$ \\
$-\Delta =t^6$

\item \tp{ {I_0^*\, II\, II}   } \quad $y^2 =x^3+t^2(2t+1)x+t^3(t^2+2t)$ \\
$-\Delta =t^9+2t^6=t^6(t^3+2)$

\item \tp{ {I_0^*\, III\, III}  } \quad $y^2 =x^3+t^2(2t+1)x$ \\
$-\Delta =t^9+2t^6=t^6(t^3+2)$

\item \tp{ {I_0^*\, II\, III}  } \quad $y^2 =x^3+t^2(2t+1)x+t^3(t+2)$ \\
$-\Delta =t^9+2t^6=t^6(t^3+2)$

\item \tp{ {I_0^*\, IV\, 1}  } \quad $y^2 =x^3+tx^2+t^2x+t^3(t+2)$ \\
$-\Delta =-t^7+t^6=-t^6(t-1)$

\item \tp{ {I_0^*\, II\, 2}  } \quad $y^2=x^3+tx^2+t^3(t^2+t+1)$\\
$-\Delta =t^8+t^7+t^6=t^6(t-1)^2$

\item \tp{ {I_0^*\, II\, 1^2}  } \quad $y^2=x^3+tx^2+t^3(t^2-1)$\\
$-\Delta =t^8-t^6=t^6(t-1)(t+1)$

\item \tp{ {I_0^*\, III\, 3}  } \quad $y^2=x^3+tx^2+t^6+t^3$\\
$-\Delta =t^9+t^6 = t^6(t+1)^3$

\item \tp{ {I_0^*\, III\,  2 \, 1 }   } \quad $y^2 =x^3+tx^2i+t^2(2t+1)x$ \\
$-\Delta =t^9+2t^8+2t^7+t^6=t^6(t-1)^2(t+1)$

\item \tp{ {I_0^*\, III\, 1^3}   } \quad $y^2 =x^3+tx^2+t^2(t+1)x+2t^3$ \\
$-\Delta =2t^9+t^8+2t^7+t^6=-t^6(t-1)\left(t-\frac{1}{i}\right)\left(t+\frac{1}{i}\right)$

\item \tp{ {I_0^*\, II\, 3}   } \quad $y^2 =x^3+tx^2+t^2(t+1)x+t^3(t^2+2t+2)$ \\
$-\Delta =2t^9+t^6=-t^6(t-1)^3$

\item \tp{ {I_0^*\, II\,  2 \, 1 }  } \quad $y^2 =x^3+tx^2+t^2(2t+1)x+t^3(2t^2+2t+2)$ \\
$-\Delta =t^9+2t^8+2t^7+t^6=t^6(t-1)^2(t+1)$

\item \tp{ {I_0^*\, II\, 1^3}   } \quad $y^2 =x^3+tx^2i+t^2(t+1)x+t^3(t^2i-ti)$ \\
$-\Delta =2t^9+t^8+2t^7+t^6 =-t^6(t-1)(t-i)(t+i)$

\end{enumerate}

\subsubsection{Case 5A}
$$\delta=7,\quad r=5\quad \mbox{Type \tp{I_1^*}} $$
\begin{equation}
y^2=x^3+tc_1x^2+t^4c_0x+t^4c_2, \quad t\nmid c_1, t\nmid c_2
\end{equation}
$$-\Delta=t^2c_1^2(t^5c_1c_2-t^8c_0^2)+t^{12}c_0^3$$

\begin{enumerate}
\itemsep=.5ex
\item \tp{ {I_1^*\, 5}  } \quad does not exist,~\ref{sigma-r}

\item \tp{ {I_1^*\,  4 \, 1 }  } \quad does not exist,~\ref{algman}

\item \tp{ {I_1^*\,  3 \, 2 }  } \quad does not exist,~\ref{lattice}

\item \tp{ {I_1^*\,  3 \, 1^2}   } \quad $y^2 =x^3+tx^2+2t^4x+t^4(t^2+1)$ \\
$-\Delta =t^2(t^5(t^2+1)+2t^8)+2t^{12})=t^7(t^2+1)(2t^3+1)$

\item \tp{ {I_1^*\,  2 \, 1^3}  } \quad $y^2 =x^3+t(t+1)x^2+it^4x+t^4((1+i)t+i)$ \\
$-\Delta =t^2(t+1)^2(t^5(t+1)((1+i)t+i)+t^8)-t^{12}i)=t^7(2+i)(2t^3+it^2+it+t+2+i)(t^2+t+1)$

\item \tp{ {I_1^*\, 2^2 \, 1 }  } \quad $y^2 =x^3+t(t+1)x^2+t^4\alpha x+t^4((2\alpha^3+2\alpha^2)t+2\alpha^3)$ \\
$-\Delta=t^2(t+1)^2(t^5(t+1)((2\alpha^3+2\alpha^2)+2\alpha^3)+2t^8\alpha^2)+t^{12}\alpha^3=t^7\alpha^2(2\alpha t^4+2\alpha t^3+2t^3+2t^5+2\alpha t+2t+2\alpha+\alpha t^5)$, where  $\alpha^4+\alpha+1=0$

\item \tp{ {I_1^*\, 1^5}   } \quad $y^2 =x^3+t(t+1)x^2+t^4x+t^4(t^2+2t+1)$ \\
$-\Delta =t^2(t+1)^2(t^5(t+1)(t^2+2t+1)+2t^8)+t^{12})=t^7(t^5+t^2+2t+1)$

\item \tp{ {I_1^*\, II\, 1^2}  } \quad $y^2 =x^3+t(t+2)x^2+t^4(t^2+1)$ \\
$-\Delta =t^7(t+2)^3(t^2+1)$

\item \tp{ {I_1^*\, II\, 2}  } \quad $y^2 =x^3+t(t+2)x^2+t^4(2t^2+t+2)$ \\
$-\Delta =t^7(t+2)^3(2t^2+t+2)$

\item \tp{ {I_1^*\, II\, 1}   } \quad $y^2 =x^3+t(t+2)x^2+t^4(2t^2+1)$ \\
$-\Delta =t^7(t+2)^3(2t^2+1)$

\item \tp{ {I_1^*\, IV}  } \quad $y^2 =x^3+t(t+2)x^2+t^4(t^2+t+1)$ \\
$-\Delta =t^7(t+2)^3(t^2+t+1)$

\item \tp{ {I_1^*\, III\, 2}  } \quad  does not exist,~\ref{algman}
\end{enumerate}

\subsubsection{Case 5B}
$$\delta=8,\quad r=6\quad \mbox{Type \tp{I_2^*}} $$
\begin{equation}
y^2=x^3+tc_1x^2+t^4c_0x+t^5d_1, \quad t\nmid c_1, t\nmid d_1
\end{equation}
$$-\Delta= t^2c_1^2(t^6c_1d_1-t^8c_0^2)+t^{12}c_0^3$$

\begin{enumerate}
\itemsep=.5ex
 \item \tp{ {I_2^*\, 4}  } \quad does not exist,
~\ref{sigma-r}

\item \tp{ {I_2^*\,  3 \, 1 }  } \quad does not exist,~\ref{lattice}

\item$\mathrm {I_2^*\, 2^2}$  \quad $y^2 =x^3+tx^2+2t^4x+2t^5$ \\
$-\Delta =t^{12}+t^{10}+t^8=t^8(t-1)^2(t+1)^2$

\item$\mathrm {I_2^*\, 21^2}$ \quad $y^2 =x^3+tx^2+t^4x+t^5(t+2)$ \\
$-\Delta =2t^{12}+t^{10}+2t^9+t^8=t^8(t-1)^2(2t^2+t+1)$

\item \tp{ {I_2^*\, 1^4}   } \quad $y^2 =x^3+tx^2+t^4x+t^5(2t+1)$ \\
$-\Delta =2t^{12}+t^{10}+t^9+2t^8=t^8(t-1)(2t^3+t^2+1)$

\item \tp{ {I_2^*\, II}  } \quad $y^2 =x^3+t(t+2)x^2+t^5(2t+1)$ \\
$-\Delta =t^8(t+2)^3(2t+1)$

\item \tp{ {I_2^*\, II\, 1}  } \quad $y^2 =x^3+t(t+2)x^2+t^5(2t+2)$ \\
$-\Delta =t^8(t+2)^3(2t+2)$

\end{enumerate}

\subsubsection{Case 5C}
$$\delta=9,\quad r=7\quad \mbox{Type \tp{I_3^*}} $$
\begin{equation}
y^2=x^3+tc_1x^2+t^4c_0x+t^6d_0, \, t\nmid c_1, \, t\nmid d_0
\end{equation}
$$-\Delta=t^2c_1^2(t^7c_1d_0-t^8c_0^2)+t^{12}c_0^3$$

\begin{enumerate}
\itemsep=.5ex
\item \tp{ {I_3^*\, 3}  } \quad does not exist,~\ref{sigma-r}

\item \tp{ {I_3^*\,  1 \, 2 }  } \quad does not exist,~\ref{lattice}

\item \tp{ {I_3^*\, 1^3}  } \quad $y^2 =x^3+tx^2+2t^4x+2t^6$ \\
$-\Delta =t^{12}+t^{10}+t^9=t^9(t-1)(t^2+t+2)$

\item \tp{ {I_3^*\, II}  } \quad $y^2 =x^3+t(t+2)x^2+2t^6$ \\
$-\Delta =2t^9(t+2)^3$

\item \tp{ {I_3^*\, III}  } \quad does not exist,~\ref{lattice}
\end{enumerate}

\subsubsection{Case 5D}
$$\delta=10,\quad r=8\quad \mbox{Type \tp{ I_4^*}}$$
\begin{equation} y^2=x^3+tc_1x^2+t^4c_0x\quad t\nmid c_0, t\nmid
c_1
\end{equation}
$$-\Delta=-t^{10}c_1^2c_0^2+t^{12}c_0^3$$

\begin{enumerate}
\itemsep=.5ex
\item \tp{ {I_4^*\, 2}  } \quad  does not exist,~\ref{sigma-r}

\item \tp{ {I_4^*\, 1^2}  } \quad $y^2 =x^3+tx^2+t^4x$ \\
$-\Delta =t^{10}(t^2-1)=t^{10}(t+1)(t-1)$

\end{enumerate}

\subsubsection{Case 6A}

$$\delta=9,\quad r=6\quad \mbox{Type \tp{IV^*}} $$
\begin{equation}
y^2=x^3+t^2c_0x^2+t^3c_1x+t^4c_2,\quad t\nmid c_2, t\nmid c_1
\end{equation}
$$-\Delta=t^4c_0^2(t^6c_0c_2-t^6d_1^2)+t^9c_1^3$$

\begin{enumerate}
\itemsep=.5ex
\item \tp{\tp{IV^*\, 3}  } \quad $y^2 =x^3+t^2x^2+t^3x+t^4(2t^2+1)$ \\
$-\Delta =t^{12}-t^9=t^9(t-1)^3$

\item \tp{\tp{IV^*\,  2 \, 1 }   } \quad $y^2 =x^3+t^2x^2+2t^3x+2t^6+t^5+2t^4$ \\
$-\Delta =t^{12}+2t^{11}+2t^{10}+t^9=t^9(t-1)^2(t+1)$

\item \tp{\tp{IV^*\, 1^3}  } \quad $y^2 =x^3+t^2x^2+t^4x+2t^6$ \\
$-\Delta =t^{12}+t^{11}+2t^{10}+t^9=t^9(t^3+t^2-t+1)$

\item \tp{\tp{IV^*\, II}  } \quad $y^2 =x^3+t^3(2t+1)x+t^4$ \\
$-\Delta =t^9(2t+1)^3$

\item \tp{\tp{IV^*\, III}  } \quad $y^2 =x^3+t^3(t+2)x+t^4(t^2+t+1)$ \\
$-\Delta =t^9(t+2)^3$

\end{enumerate}

\subsubsection{Case 6B}
$$\delta=10,\quad r=6\quad \mbox{Type \tp{IV^*}} $$
\begin{equation}
y^2=x^3+t^2c_0x^2+t^4d_0x+ t^4 c_2 
\quad t\nmid c_2, t\nmid c_0
\end{equation}
$$-\Delta=t^4c_0^2(t^6c_0c_2-t^8d_0^2)+t^{12}d_0^3)$$

\begin{enumerate}
\itemsep=.5ex
\item \tp{ {IV^*\, 1^2}  } \quad $y^2 =x^3+tx^2 +t^4(2t^2+1)$\\
$-\Delta =t^{12}+2t^{10}=t^{10}(t+1)(t-1)$

\item \tp{ {IV^*\, 2}  } \quad $y^2 =x^3+t^2c_0x^2 +t^4(\alpha t^2+\alpha t + \alpha)$ \\
$-\Delta =t^{10}(c_0^3\alpha (t-1)^2)$

\end{enumerate}

\subsubsection{Case 6C}

$$\delta=12,\quad r=6\quad \mbox{Type \tp{IV^*}} $$
\begin{equation}
y^2=x^3+t^4d_0x+t^4c_2,\quad t\nmid c_2, t\nmid d_0
\end{equation}
$$-\Delta=t^{12}d_0^3$$

\subsubsection{Case 7}
$$\delta=9,\quad r=7\quad \mbox{Type \tp{III^*}} $$
\begin{equation}
y^2=x^3+t^2c_0x^2+t^3c_1x+t^5d_1,\quad t\nmid c_1
\end{equation}
$$-\Delta=t^4c_0^2(t^7c_0d_1-t^6c_1^2)+t^9c_1^3$$

\begin{enumerate}
\itemsep=.5ex
\item \tp{\tp{III^*\, 3}  } \quad does not exist,~\ref{sigma-r}

\item \tp{\tp{III^*\,  2 \, 1}  } \quad $y^2 =x^3+t^2x^2+t^3(t+1)x+t^6+t^5$ \\
$-\Delta =t^{12}+2t^{11}+2t^{10}+t^9=t^9(t+1)(t^2+t+1)$

\item \tp{\tp{III^*\, 1^3}  } \quad $y^2 =x^3+t^2x^2+t^3(t+1)x+2t^6$ \\
$-\Delta =2t^{12}+t^{11}+2t^{10}+t^9=t^9(2t+1)(t^2+1)$

\item \tp{\tp{III^*\, II}   } \quad $y^2 =x^3+t^3(2t+1)x+t^5(2+t)$\\
$-\Delta =t^9(2t+1)^3$

\item \tp{\tp{III^*\, III}  } \quad $y^2 =x^3+t^3(2t+1)x$ \\
$-\Delta =t^9(2t+1)^3$
\end{enumerate}

\subsubsection{Case 8A}
$$\delta=11,\quad r=8\quad \mbox{Type \tp{II^*}} $$
\begin{equation}
y^2=x^3+t^2c_0x^2+t^4d_0x+t^5d_1,\quad t\nmid d_1, t\nmid c_0
\end{equation}
$$-\Delta=t^4c_0^2(t^7c_0d_1-t^8d_0^2)+t^{12}d_0^3$$

\begin{enumerate}
\item \tp{\tp{II_2^*\, 1}  } \quad
$y^2 =t^2x^2+t^4x+t^5$ \\
$-\Delta =t^{11}(t+1)$

\end{enumerate}

\subsubsection{Case 8B}
$$\delta=12,\quad r=8\quad \mbox{Type \tp{II^*}} $$
\begin{equation}
y^2=x^3+t^4d_0x+t^5d_1,\quad t\nmid d_1, \nmid d_0
\end{equation}
$$-\Delta=t^{12}d_0^3$$

\section{Conclusion and Summary of Results}

 \subsection{General Results}
 Of the $372$ configurations checked above, having either multiplicative or additive
 singularities, $267$ configurations were found to exist in
 characteristic three, and $105$ were found not to exist. $227$
 configurations involving at least one additive singularity exist,
 while $68$ involving at least one additive singularity did not.  Of
 those configurations which are purely multiplicative, $40$ exist, while
 $37$ do not.

\subsection{Multiplicative Results in Comparison with other Characteristics}
\label{subsec:multresult}

As mentioned earlier, in different characteristics the additive
singularities have such different properties that their configurations
are not easily compared from characteristic to characteristic.
However, this is not a problem with purely multiplicative fibres, so
we compare them here.

There are $6$ configurations that exist in characteristic zero,
but not in characteristic three.  But every configuration that
exists in characteristic three also exists in
characteristic zero. There are $6$ configurations which exist in
characteristic two but not characteristic three, and there are $10$
configurations that exist in characteristic three but not
characteristic two.

The the following table summarizes the results.  A \checkmark
indicates that a rational elliptic surface with the indicated
fibre distribution exists in the given case and characteristic,
whereas an $X$ indicates that no such surface exists.

\begin{longtable}[c]{|lccc|}
\caption{Existence of multiplicative fibre types}
      \\
      \hline

Partition   & Characteristic 0 & Characteristic 2  & Characteristic 3 \\
\hline
\hline
\endfirsthead
        \hline
        \multicolumn{4}{|l|}{\small\slshape continued from previous page}\\
        \hline
Partition   & Characteristic 0 & Characteristic 2  & Characteristic 3 \\
    \hline
\endhead
        \hline
        \multicolumn{4}{|r|}{\small\slshape continued on next  page}\\
        \hline            \endfoot
\hline
 \endlastfoot

$1^{12}$ & \checkmark & \checkmark & \checkmark \\
\hline $2 \,1^{10}$
& \checkmark & \checkmark & \checkmark \\
\hline $3 \,1^9$ &
\checkmark & \checkmark & \checkmark \\
\hline $2^2\, 1^8$ & \checkmark
& \checkmark & \checkmark \\
\hline $4 \,1^8$ & \checkmark &
\checkmark & \checkmark \\
\hline $3\, 2\, 1^7$ & \checkmark &
\checkmark & \checkmark \\
\hline $5\, 1^7$ & \checkmark & \checkmark
& \checkmark \\
\hline $2^3\, 1^6$ & \checkmark & \checkmark &
\checkmark \\
\hline $4 \,2\, 1^6$ & \checkmark & \checkmark &
\checkmark \\
\hline $3^2 \,1^6$ & \checkmark & \checkmark &
\checkmark \\
\hline $6 \,1^6$ & \checkmark & \checkmark & \checkmark
\\
\hline $3 \,2^2 \,1^5$ & \checkmark & \checkmark & \checkmark \\
\hline $5 \,2\, 1^5$ & \checkmark & \checkmark & \checkmark \\
\hline
$4 \,3\, 1^5$ & \checkmark & \checkmark & \checkmark \\
\hline $7
\,1^5$ & \checkmark & \checkmark & \checkmark \\
\hline $2^4 \,1^4$ &
\checkmark & \checkmark & \checkmark \\
\hline $4 \,2^2 \,1^4$ &
\checkmark & \checkmark & \checkmark \\
\hline $3^2 \,2\, 1^4$ &
\checkmark & \checkmark & \checkmark \\
\hline $6 \,2 \,1^4$ &
\checkmark & \checkmark & \checkmark \\
\hline $5 \,3\, 1^4$ &
\checkmark & \checkmark & \checkmark \\
\hline $4^2 \,1^4$ &
\checkmark & \checkmark & \checkmark \\
\hline $8 \,1^4$ & \checkmark
& \checkmark & \checkmark \\
\hline $3 \,2^3\, 1^3$ & \checkmark &
\checkmark & \checkmark \\
\hline $5 \,2^2 \,1^3$ & \checkmark &
\checkmark & \checkmark \\
\hline $4\, 3\, 2\, 1^3$ & \checkmark &
\checkmark & \checkmark \\
\hline $7 \,2 \,1^3$ & \checkmark &
\checkmark & \checkmark \\
\hline $3^3\, 1^3$ & \checkmark &
\checkmark &$X$ \\
\hline $6 \,3\, 1^3$ & \checkmark & \checkmark &$X$ \\
\hline $5 \,4 \,1^3$ & \checkmark & \checkmark & \checkmark \\
\hline
$9 \,1^3$ & \checkmark & \checkmark &$X$ \\
\hline $2^5 \,1^2$ &
\checkmark &$X$ & \checkmark \\
\hline $4 \,2^3 \,1^2$ & \checkmark &$X$ &
\checkmark \\
\hline $3^2 \,2^2\, 1^2$ & \checkmark &\checkmark &
\checkmark \\
\hline $6 \,2^2\, 1^2$ & \checkmark &$X$ & \checkmark \\
\hline $5 \,3 \,2\, 1^2$ & \checkmark & \checkmark & \checkmark \\
\hline $4^2 \,2 \,1^2$ & \checkmark &$X$ & \checkmark \\
\hline $8 \,2
\,1^2$ & \checkmark &$X$ & \checkmark \\
\hline $4 \,3^2\, 1^2$ &$X$ &$X$ &$X$ 
\\
\hline $7 \,3 \,1^2$ &$X$ &$X$ &$X$ \\
\hline $6 \,4 \,1^2$ &$X$ &$X$ &$X$ \\
\hline $5^2\, 1^2$ & \checkmark & \checkmark & \checkmark \\
\hline
$10 \,1^2$ &$X$ &$X$ &$X$ \\
\hline $3 \,2^4\, 1$ & \checkmark &$X$ &
\checkmark \\
\hline $5 \,2^3\, 1$ &$X$ &$X$ &$X$ \\
\hline $4\, 3 \,2^2
\,1$ & \checkmark &$X$ & \checkmark \\
\hline $7 \,2^2\, 1$ &$X$ &$X$ &$X$ \\
\hline $3^3\, 2 \,1$ & \checkmark & \checkmark &$X$ \\
\hline $6 3 \,2\,
1$ & \checkmark &$X$ &$X$ \\
\hline $5 \,4\, 2 \,1$ &$X$ &$X$ &$X$ \\
\hline $9
\,2 \,1$ &$X$ &$X$ &$X$ \\
\hline $5 \,3^2\, 1$ &$X$ &$X$ &$X$ \\
\hline $4^2
\,3\, 1$ &$X$ &$X$ &$X$ \\
\hline $8\, 3 \,1$ &$X$ &$X$ &$X$ \\
\hline $7\, 4\, 1$
&$X$ &$X$ &$X$ \\
\hline $6 \,5\, 1$ &$X$ &$X$ &$X$ \\
\hline $11 \,1$ &$X$ &$X$ &$X$ \\
\hline $2^6$ & \checkmark &$X$ & \checkmark \\
\hline $4 \,2^4$ &
\checkmark &$X$ & \checkmark \\
\hline $3^2 \,2^3$ &$X$ &$X$ &$X$ \\
\hline $6
\,2^3$ &$X$ &$X$ &$X$ \\
\hline $5 \,3\, 2^2$ &$X$ &$X$ &$X$ \\
\hline $4^2 \,2^2$
& \checkmark &$X$ & \checkmark \\
\hline $8 \,2^2$ &$X$ &$X$ &$X$ \\
\hline
$4\, 3^2\, 2$ &$X$ &$X$ &$X$ \\
\hline $7 \,3 \,2$ &$X$ &$X$ &$X$ \\
\hline $6\, 4
\,2$ &$X$ &$X$ &$X$ \\
\hline $5^2\, 2$ &$X$ &$X$ &$X$ \\
\hline $10 \,2 $ &$X$ &$X$ 
&$X$ \\
\hline $3^4$ & \checkmark & \checkmark &$X$ \\
\hline $6 \,3^2$ &$X$ 
&$X$ &$X$ \\
\hline $5 \,4 \,3$ &$X$ &$X$ &$X$ \\
\hline $9 \,3$ &$X$ &$X$ &$X$ \\
\hline $4^3$ &$X$ &$X$ &$X$ \\
\hline $8 \,4$ &$X$ &$X$ &$X$ \\
\hline $7 \,5$ &$X$ 
&$X$ &$X$ \\
\hline $6^2$ &$X$ &$X$ &$X$ \\
\hline $12$ &$X$ &$X$ &$X$ \\
\hline

\end{longtable}

\typeout{LaTeX Warning: Label(s) may have changed. Rerun}


\newpage




\begin{thebibliography}{www}

\bibitem{dynkin} E. B. Dynkin, \emph{Semisimple subalgebras of
    semisimple Lie algebras}, Mat. Sbornik N.S., \textbf{Vol. 30(72)}
  1952, 349-462.

\bibitem{Lang-I} W. E. Lang, \emph{Configurations of singular fibres
on rational elliptic surfaces in characteristic two},
Communications in Algebra, \textbf{Vol. 28(12)} 2000, 5813--5836.

\bibitem{Lang-II} W. E. Lang, \emph{Extremal rational elliptic surfaces in
characteristic $p$. II:Surfaces with three or fewer singular
fibres}, Arkiv for matematik, \textbf{Vol. 32} 1994, 423--438.

\bibitem{miranda-singfibres} R. Miranda, \emph{Persson's list of
singular fibers for a rational elliptic surface}, Mathematische
Zeitschrift, \textbf{Vol. 205} 1990, 191--211.

\bibitem{mirper} R. Miranda and U. Persson, \emph{On extremal rational
    elliptic surfaces}, Mathematische Zeitschrift, \textbf{Vol. 193}
1986, 537--558.

\bibitem{oguiso-shioda} K. Oguiso and T. Shioda, \emph{The
    Mordell-Weil lattice of a rational elliptic surface}, Comment.
  Math. Univ. St. Paul, \textbf{Vol. 40(1)} 1990, 83--99.

\bibitem{persson-fibreconfig} U. Persson, \emph{Configurations of
{K}odaira fibers on rational elliptic surfaces}, Mathematische
Zeitschrift, \textbf{Vol. 205(1)} 1990, 1--47.

\bibitem{shioda} T. Shioda, \emph{On the Mordell-Weil lattices}, Comment. Math.
  Univ. St. Paul, \textbf{Vol. 39(2)} 1990, 211--240.

\bibitem{Tate} J. H. Silverman, \emph{Advanced topics in the
arithmetic of elliptic curves}, Graduate Texts in Mathematics,
\textbf{Vol. 151} 1994.



\end{thebibliography}
\end{document}